\documentclass[11pt]{article}

\usepackage[T1]{fontenc}
\usepackage[utf8]{inputenc}
\usepackage{lmodern}
\usepackage[a4paper,margin=1in]{geometry}
\usepackage{amsmath,amssymb,amsthm}
\usepackage{mathtools}
\usepackage{enumitem}
\usepackage[numbers,sort&compress]{natbib}
\usepackage[colorlinks=true,linkcolor=blue,citecolor=blue,urlcolor=blue]{hyperref}
\usepackage[nameinlink,noabbrev]{cleveref}

\newcommand{\C}{\mathbb{C}}

\newtheorem{proposition}{Proposition}[section]
\newtheorem{theorem}{Theorem}[section]
\newtheorem{lemma}{Lemma}[section]

\title{Universality of kernels on Riemannian symmetric spaces}

\author{%
Salem Said \\
\small Laboratoire Jean-Kuntzmann (CNRS), Universit\'e Grenoble-Alpes, Saint Martin d'H\`eres, France
\and
Natha\"el Da Costa\\
\small T\"ubingen AI Center, University of T\"ubingen, T\"ubingen, Germany
\and
Franziskus Steinert \\ Cyrus Mostajeran\\
\small Division of Mathematical Sciences, Nanyang Technological University, Singapore
}

\date{}

\begin{document}

\maketitle

\begin{abstract}
We investigate universality properties of continuous, positive-definite invariant kernels on Riemannian symmetric spaces, providing a unified harmonic-analytic characterization across compact and non-compact settings. In the compact symmetric case, we prove that a continuous, positive-definite invariant kernel is $C$-universal if and only if all of its spherical coefficients are strictly positive, a sharpening of classical Bochner-type results. This characterization is extended to compact homogeneous spaces, where universality is shown to be equivalent to the strict positive-definiteness of coefficient matrices arising from a representation-theoretic expansion. In contrast, for non-compact symmetric spaces, we establish that any continuous, positive-definite invariant kernel that is $C_0$ and integrable is automatically $C_0$-universal. Our analysis essentially relies on spectral decompositions using spherical functions and group Fourier transforms, in order to provide alternative, harmonic-analytic formulations of universality conditions. Several examples (including kernels on spheres and compact Lie groups, as well as integrable kernels on hyperbolic spaces and symmetric cones) both illustrate the theory and demonstrate practical criteria for constructing universal kernels. 
\end{abstract}


\section{Introduction}
Positive-definite functions are fundamental objects in both harmonic analysis and machine learning. In the former field, a positive-definite function on a group uniquely defines a unitary representation of this group. In the latter, positive-definite functions are at the heart of kernel methods, an immensely popular class of learning methods. Both~fields rely on the same mathematical construction, which associates to a positive-definite function (or kernel) a reproducing-kernel Hilbert space of bounded continuous functions~\cite{paulsen}.

In machine learning, a fundamental question is to determine, given a kernel $\mathcal{K}$ on some locally-compact topological space $X$, whether this kernel is universal or not. Universality means the reproducing-kernel Hilbert space $H$ associated with $\mathcal{K}$ is sufficiently rich, and can therefore be used to approximate ``arbitrary functions''~\cite{steinwart2001influence}\cite{christmann_universal_2010}\cite{sriperumbudur2011universality}\cite{carmeli}. 

Of course, there are several versions of the general concept of universality. When $X$ is compact, $\mathcal{K}$ is called $C$-universal if $H$ is dense in $C(X)$ (the space of continuous functions on $X$, which has the uniform convergence topology). When $X$ is non-compact, $\mathcal{K}$ is called $C_0$-universal if $H$ is dense in $C_0(X)$ (the space of continuous functions that vanish at infinity, also equipped with the uniform convergence topology). 

In~\cite{sriperumbudur2011universality}, a necessary and sufficient condition was given, for both of these kinds of universality. A kernel $\mathcal{K}$ on $X$ is~called integrally strictly positive-definite, if it satisfies
\begin{equation} \label{eq:integralpd}
 \int_X\int_X\,\mathcal{K}(x,y)\mu(dx)\mu(dy) > 0
\end{equation}
for all non-zero finite signed measures $\mu$ on $X$. According to~\cite{sriperumbudur2011universality} (Proposition 4), $\mathcal{K}$ is $C$-universal (when $X$ is compact), or $C_0$-universal (when $X$ is locally-compact, but noncompact) if and only if $\mathcal{K}$ is integrally strictly positive-definite. It~should be emphasised that this result only applies to real kernels ($\mathcal{K}:X\times X \rightarrow \mathbb{R}$) and to real function spaces (in other words, $H$, $C(X) = C(X,\mathbb{R})$ or $C_0(X) = C_0(X,\mathbb{R})$ are spaces of real functions). 

In case $X$ is an abelian group (say $X = \mathbb{R}^d$) and $\mathcal{K}$ is a translation-invariant kernel, so $\mathcal{K}(x,y) = f_\mathcal{K}(x-y)$ where $f_\mathcal{K}$ is a positive-definite function, Condition (\ref{eq:integralpd}) can be rewritten in the Fourier domain~\cite{carmeli}. Indeed, Bochner's theorem says that $f_\mathcal{K}$ is the inverse Fourier transform of a finite positive measure $\Lambda$, and Condition (\ref{eq:integralpd}) turns out to be equivalent to $\Lambda$ having full support\,: 
$\mathrm{supp}\,\Lambda = \hat{X}$ where $\hat{X}$ is the dual of $X$.

The general aim of the present work is to extend this kind of result to the case where $X$ is a Riemannian symmetric space, either of compact or non-compact type. In doing so, we both significantly improve and extend preliminary investigations pursued in \cite{Steinart2026} concerned with establishing sufficient conditions for universality in the non-compact case. The case of compact homogeneous spaces, which are not necessarily symmetric, is also considered. This is motivated by the application of kernel methods to learning from non-Euclidean data. Specifically, many recent applications, such as medical image analysis, brain-computer interface analysis, radar signal processing, or graph representation learning, feature data that belong to non-Euclidean spaces \cite{congedo}\cite{pennec_book}\cite{radar}\cite{graphembedding}, which (in all of the applications just mentioned) turn out to be Riemannian symmetric spaces. 

Roughly, the idea will be to use the machinery of harmonic analysis on symmetric spaces, in order to represent a kernel $\mathcal{K}$ on a symmetric space $X$ in the Fourier domain (or to obtain a ``spectral decomposition" of such a kernel). This is a first step to rewriting Condition (\ref{eq:integralpd}) in an alternative form which yields new characterizations of universality. The above-said machinery essentially relies on spherical transforms and group Fourier transforms~\cite{helgason2022groups}\cite{kirillov}.  

In the following, $X = G/K$ denotes a Riemannian symmetric space, either of compact type (Section \ref{sec:compact}) or of non-compact type (Section \ref{sec:noncompact}). Thus, the isometry group $G$ is a semisimple Lie group and $K$ a compact subgroup of $G$. An invariant kernel $\mathcal{K}: X \times X \rightarrow \mathbb{C}$ is uniquely determined by the function $f_\mathcal{K}:X \rightarrow \mathbb{C}$, given by $f_\mathcal{K}(x) = \mathcal{K}(x,K)$ (note that $X$ is understood as a space of left cosets $x = gK$ for $g \in G$). Indeed,
\begin{equation} \label{eq:basicidentity}
    \mathcal{K}(g_1K,g_2K) = f_\mathcal{K}(g^{-1}_2g^{\phantom{-1}}_1\!\!\!\!K) \hspace{0.5cm}
\text{for any $g_1, g_2 \in G$}
\end{equation}
The case of symmetric spaces of compact type is considered in Section \ref{sec:compact}. To begin, the function $f_\mathcal{K}$ is expanded into a spherical series
\begin{equation} \label{eq:spherical_series_intro}
    f_\mathcal{K}(x) = c_X\sum_{\lambda \in \Lambda} d_\lambda\hspace{0.02cm}\hat{f}_\mathcal{K}(\lambda)\hspace{0.03cm}\varphi_\lambda(x)
\end{equation}
where $c_X > 0$ is a normalising constant, and the sum is over all irreducible spherical representations $\lambda$ of $G$. 
Each~representation $\lambda$ has dimension $d_\lambda$ and is uniquely determined by the (positive-definite) spherical function $\varphi_\lambda\hspace{0.03cm}$. The~coefficients $\hat{f}_\mathcal{K}(\lambda)$ are the spherical coefficients of the function $f_\mathcal{K}\hspace{0.03cm}$. The spectral decomposition of the kernel $\mathcal{K}$ then follows from the basic identity (\ref{eq:basicidentity}) and from unitarity of the irreducible spherical representations 
(see~(\ref{eq:compact_spectral_decompo}) and (\ref{eq:compact_spectral_decompo_bis}) in Paragraph \ref{subsec:proof_c_lemma}). This spectral decomposition is crucial to the proof of the main result of Section \ref{sec:compact}, Theorem \ref{th:compact_univ}\,:\hfill\linebreak the kernel $\mathcal{K}$ is $C$-universal if and only if all the spherical coefficients $\hat{f}_\mathcal{K}(\lambda)$ are strictly positive. This result is thus analogous to the known result in Euclidean spaces.

Section \ref{sec:compact_hom} extends the approach of Section \ref{sec:compact} to compact homogeneous spaces, which are not necessarily symmetric. The main result is Theorem \ref{th:compact_hom_univ}, which is a direct generalisation of Theorem \ref{th:compact_univ}.

Section \ref{sec:noncompact} pursues the same template for symmetric spaces of non-compact type. It turns out that this requires new methods, and also leads to a different conclusion. In some sense, guaranteeing universality in this case requires even less than in the Euclidean and compact cases. This section assumes that the kernel $\mathcal{K}$ is integrable, which means the function $f_\mathcal{K}$ is integrable with respect to the Riemannian volume measure of $X$. This assumption guarantees $f_\mathcal{K}$ can be expressed as an inverse spherical transform
\begin{equation} \label{eq:L1god_intro}
  f_\mathcal{K}(x) = c_X\int_{\mathfrak{a}^*}\,\hat{f}_\mathcal{K}(\lambda)\,\varphi_\lambda(x)\hspace{0.02cm}\sigma(d\lambda)    
\end{equation}
which only involves spherical functions $\varphi_\lambda$ which correspond to irreducible unitary representations from the so-called spherical principal series of representations of $G$ (see Paragraph \ref{subsec:nctools}). Here, the parameter $\lambda$ belongs to the dual $\mathfrak{a}^*$ of~the~(real) Lie algebra $\mathfrak{a}$ of the abelian Lie group $A$ appearing in the Iwasawa decomposition $G = KAN$~\cite{helgason2022groups}. Moreover, $\sigma(d\lambda)$ is a spectral measure involving the Harish-Chandra function (details are provided in Paragraph \ref{subsec:noncompact}).    
The expression (\ref{eq:L1god_intro}) is a result of the $L^1$-Godement theorem, which moreover states that $\hat{f}_\mathcal{K}(\lambda) \geq 0$ for all $\lambda \in \mathfrak{a}^*$~\cite{dacosta2024invariantkernelsriemanniansymmetric}. The function $\hat{f}_\mathcal{K}(\lambda)$ is the spherical transform of $f_\mathcal{K}\hspace{0.03cm}$. As explained in~\cite{dacosta2024invariantkernelsriemanniansymmetric}, the $L^1$-Godement theorem is a special case of the so-called Bochner-Godement theorem, which describes all positive-definite functions (not just integrable ones).

Pursuing the approach of Section \ref{sec:compact}, (\ref{eq:L1god_intro}) is then used to obtain a spectral decomposition of the kernel $\mathcal{K}$ (see (\ref{eq:kernel_IFT3})), rewrite Condition (\ref{eq:integralpd}) (see Paragraph \ref{subsec:noncomp-lem_prooof}), and arrive at a final conclusion in the form of Theorem \ref{th:ncc0}\,: the kernel $\mathcal{K}$ is $C_0$-universal whenever it is integrable. This is quite different from the conclusion made in Theorem \ref{th:compact_univ}, and this difference is above all due to the inherent analyticity of the spherical transform $\hat{f}_\mathcal{K}(\lambda)\hspace{0.03cm}$. Indeed (as recalled in Lemma \ref{lem:holomorphic}, based on~\cite{helgason2022groups}), $\hat{f}_\mathcal{K}(\lambda)\hspace{0.03cm}$ extends to a holomorphic function on a certain tube domain in the complexification of $\mathfrak{a}^*$. 

Sections \ref{sec:compact}, \ref{sec:compact_hom} and \ref{sec:noncompact} all operate under the assumption that the kernel $\mathcal{K}$ is real. In practice, there is nothing lost by using a real (rather than complex) kernel. If the reproducing-kernel Hilbert space $H$ associated to a real kernel $\mathcal{K}$ is a dense subspace of $C(X)$ or $C_0(X)$, then the complexification of $H$ is dense in the corresponding space of complex-valued functions ($C(X,\mathbb{C})$ or $C_0(X,\mathbb{C})$).

Section \ref{sec:grouptospherical} expands on the manner in which the group Fourier transform was employed in Section \ref{sec:noncompact}. It derives explicit inversion, convolution, and Plancherel formulas for the group Fourier transform, when applied to symmetric spaces of non-compact type. While the Fourier-Helgason transform is usually considered the standard ``Fourier transform'' for functions on symmetric spaces of non-compact type, Section \ref{sec:grouptospherical} explains that it is not suitable for the purpose of the present work, and that the group Fourier transform is the more adequate choice.     

The final Section \ref{sec:examples} illustrates how the results from Sections \ref{sec:compact} and \ref{sec:noncompact} can be used to construct universal kernels on spheres, compact Lie groups, symmetric cones, and type IV symmetric spaces (the dual symmetric spaces of semisimple compact Lie groups). The examples obtained in this section showcase the ability to systematically construct easy-to-evaluate (in fact, closed-form) universal kernels. This is complimentary to existing work in the machine learning literature, which primarily focuses on  numerical evaluation of heat or Matérn kernels~\cite{azangulov2023stationary}\cite{azangulov2024stationary}.  

The main theorems in the present work (\ref{th:compact_univ}, \ref{th:compact_hom_univ}, and \ref{th:ncc0}) can be thought of as Tauberian theorems. The reproducing-kernel Hilbert space associated with the kernel $\mathcal{K}$ is spanned by left translates of the function $f_\mathcal{K}$ (these are the functions $x \mapsto f_\mathcal{K}(gx)$ where $g \in G$). Universality means the span of these translates is dense in $C(X)$ or $C_0(X)$. Tauberian theorems on symmetric spaces of non-compact type have previously been considered in~\cite{sitaram2011}, but these are concerned with density of the span of the translates in $L^1(X,\mathrm{vol})$, where $\mathrm{vol}$ denotes the Riemannian volume measure.    

\section{Compact symmetric spaces} \label{sec:compact}

\subsection{Statement of the main result} \label{subsec:compstatement}
In this section, $X = G/K$ will be a Riemannian symmetric space of the compact type. The~aim is to prove the following result\,: an invariant kernel on $X$ is $C$-universal if and only if its spherical coefficients are all strictly positive.

A kernel on a set $X$ is a Hermitian map $\mathcal{K}:X\times X \rightarrow \mathbb{C}$. The kernel $\mathcal{K}$ is said to be positive definite if for all $N\in\mathbb{N}$ and $x_1,\dots,x_N\in X$,
the Gram matrix $\big[\mathcal{K}(x_i,x_j)\big]_{i,j}$ is Hermitian positive semidefinite.
It is said to be invariant if $\mathcal{K}(g\hspace{0.02cm}x,g\hspace{0.02cm}y) = \mathcal{K}(x,y)$ for any $g \in G$ ($g\hspace{0.02cm} x$ is the action of $g \in G$ on $x \in G/K$ by multiplication from the left). To an invariant kernel $\mathcal{K}$, it~is~useful to associate the function $f_\mathcal{K}(x) = \mathcal{K}(x,K)$. This is a $K$-invariant function, $f_\mathcal{K}(k\hspace{0.02cm}x) = f_\mathcal{K}(x)$ for any $k \in K$.
The spherical coefficients of the kernel $\mathcal{K}$ appear in the expansion of $f_\mathcal{K}$ in terms of spherical functions on $X$. 

Here, it is helpful to recall some background about this expansion.  
Let $\Lambda$ denote the set of equivalence classes of spherical irreducible representations of $G$ (this is a countable set~\cite{helgason2022groups}). Each $\lambda \in \Lambda$ contains a unitary representation $U_\lambda$ in a Hilbert space $V_\lambda$ of dimension $d_\lambda <\infty$, and admits a unique (up~to~sign) unit vector $h_\lambda$ such that $U_\lambda(k)\hspace{0.02cm}h_\lambda = h_\lambda$ for any $k \in K$. The corresponding spherical function is
\begin{equation} \label{eq:compact_spherical}
          \varphi_\lambda(gK) = \langle h_\lambda\hspace{0.02cm},U_\lambda(g)h_\lambda\rangle
       \hspace{0.5cm}\text{(scalar product in $V_\lambda$)}
\end{equation}
For any continuous $K$-invariant function $f:X \rightarrow \mathbb{C}$, one has the spherical coefficients
\begin{equation} \label{eq:spherical_coeff}
    \hat{f}(\lambda) = \int_X\,f(x)\hspace{0.02cm}\varphi_\lambda^*(x)\hspace{0.03cm}\mathrm{vol}(dx)
\end{equation}
where $^*$ denotes the complex conjugate, and $\mathrm{vol}$ the volume measure of $X$. These yield a spherical series, 
\begin{equation} \label{eq:spherical_series}
    f(x) = c_X\,\sum_{\lambda \in \Lambda} d_\lambda\hspace{0.02cm}\hat{f}(\lambda)\hspace{0.03cm}\varphi_\lambda(x)
\end{equation}
where $c_X$ is a normalising constant. The spherical series (\ref{eq:spherical_series}) is convergent in $L^2(X,\mathrm{vol})$ and it is moreover summable (using suitable resummation methods) in the uniform norm of $C(X,\mathbb{C})$~\cite{helgason2022groups}. 
\begin{theorem} \label{th:compact_univ}
    A continuous, positive-definite invariant kernel $\mathcal{K}:X\times X \rightarrow \mathbb{R}$ is $C$-universal if and only if $\hat{f}_\mathcal{K}(\lambda) > 0$ for each $\lambda \in \Lambda$ (where the spherical coefficients $\hat{f}_\mathcal{K}(\lambda)$ are given as in (\ref{eq:spherical_coeff})). 
\end{theorem}
This theorem is a sharpening of Bochner's theorem for compact symmetric spaces (see~\cite{yaglom_second-order_1961}), which
states that a continuous invariant kernel $\mathcal{K}:X\times X \rightarrow \mathbb{C}$ is positive-definite if and only if $\hat{f}_\mathcal{K}(\lambda) \geq 0$ for each $\lambda \in \Lambda$. 

In order to prove Theorem \ref{th:compact_univ}, Condition (\ref{eq:integralpd}) will be rewritten in an alternative form, involving the spherical representations $U_\lambda\hspace{0.03cm}$. Let $\mu$ be a finite signed measure on $X$ and consider for each $\lambda \in \Lambda$, the vector $\hat{\mu}(\lambda) \in V_\lambda$,
\begin{equation} \label{eq:compact_fourier}
    \hat{\mu}(\lambda) = \int_X\,(U_\lambda\hspace{0.03cm}h_\lambda)(x)\hspace{0.03cm}\mu(dx)
\end{equation}
where the notation $(U_\lambda\hspace{0.03cm}h_\lambda)(x)$ has the following meaning\,: since $U_\lambda(k)\hspace{0.02cm}h_\lambda = h_\lambda$ for any $k \in K$, $U_\lambda(g)\hspace{0.02cm}h_\lambda$ is constant on any coset $x = gK$, and therefore defines a function on $X$.

The proof of Theorem \ref{th:compact_univ} relies on the following lemma.
\begin{lemma} \label{lem:integralpd_four_comp}
    Let $\mathcal{K}$ be an invariant kernel on $X$ and $f_\mathcal{K}$ the associated function, as above. Then, 
    \begin{equation} \label{eq:integralpd_four_comp}
       \text{LHS of (\ref{eq:integralpd})} = c_X\,\sum_{\lambda \in \Lambda}\,d_\lambda\hspace{0.03cm}\hat{f}_\mathcal{K}(\lambda)\hspace{0.03cm}\Vert \hat{\mu}(\lambda)\Vert^2
    \end{equation}
\end{lemma}
Moreover, the proof uses the fact that $\mu = 0$ (the zero measure) if and only if $\hat{\mu}(\lambda) = 0$ for all $\lambda \in \Lambda$. This is equivalent to the statement that the functions
$$
\varphi^h_\lambda(x) = \langle h,(U_\lambda\hspace{0.03cm}h_\lambda)(x)\rangle \hspace{0.5cm} \text{ where $h \in V_\lambda$}
$$
span a dense subset of $C(X,\mathbb{C})$ (with respect to the uniform norm). In turn, this is from~\cite{helgason2022groups} (Theorem 4.3, Page 538). 

\subsection{Proof of the main result}
Admit for now that Lemma \ref{lem:integralpd_four_comp} is true. Here is then a proof of Theorem \ref{th:compact_univ}.

Assume $\hat{f}_\mathcal{K}(\lambda) > 0$ for all $\lambda \in \Lambda$. If $\mu$ is any non-zero finite signed measure, then $\hat{\mu}(\lambda) \neq 0$ for some $\lambda \in \Lambda$ 
(by~the~fact mentioned after Lemma \ref{lem:integralpd_four_comp}). The corresponding term in the sum on the right-hand side of (\ref{eq:integralpd_four_comp}) is therefore strictly positive, which implies that Condition (\ref{eq:integralpd}) is satisfied. In~other words, $\mathcal{K}$ is $C$-universal.

Conversely, assume $\mathcal{K}$ is $C$-universal, so Condition (\ref{eq:integralpd}) holds true. Bochner's theorem for compact symmetric spaces (mentioned after Theorem \ref{th:compact_univ})  says the spherical coefficients $\hat{f}_\mathcal{K}(\lambda)$ are either positive or zero. If there did exist $\lambda_0 \in \Lambda$ with $\hat{f}_\mathcal{K}(\lambda_0) = 0$, then it would be possible to construct a non-zero finite signed measure $\mu$ that violates Condition (\ref{eq:integralpd}).

Indeed, consider the function $m:X \rightarrow \mathbb{C}$ given as a convolution $m = \varphi^{\phantom{*}}_{\lambda_0} * \varphi^*_{\lambda_0}$ where $^*$ denotes the complex conjugate. The convolution operation is the usual group convolution~\cite{helgason2022groups}\cite{kirillov}. A very important feature of symmetric spaces is that the convolution of $K$-invariant functions is a commutative operation~\cite{helgason2022groups} (Page 292). In particular, $\varphi^{\phantom{*}}_{\lambda_0}*\varphi^*_{\lambda_0} = 
\varphi^*_{\lambda_0}*\varphi^{\phantom{*}}_{\lambda_0}$ which shows that $m$ is real, in addition to being $K$-invariant.\\[0.1cm]
\indent Consider then the finite signed measure $\mu (dx)= m(x)\hspace{0.02cm}\mathrm{vol}(dx)$. To show that this violates Condition (\ref{eq:integralpd}), it is enough to show the right-hand side of (\ref{eq:integralpd_four_comp}) is equal to zero (for this $\mu$). 
To this end, it will be proved below that
\begin{equation} \label{eq:proof_comp_cov_1}
  \hat{\mu}(\lambda) = \delta_{\lambda\lambda_0}\hspace{0.03cm}c(\lambda_0)\hspace{0.03cm}h_{\lambda_0}    
\end{equation}
where $\delta_{\lambda\lambda_0} = 1$ if $\lambda = \lambda_0$ and $\delta_{\lambda\lambda_0}= 0$ otherwise, and where $c(\lambda_0)$ is a real number. However, replacing (\ref{eq:proof_comp_cov_1}) into (\ref{eq:integralpd_four_comp}), the~result is clearly zero. 

Since $\mu$ is non-zero and violates Condition (\ref{eq:integralpd}), $\mathcal{K}$ cannot be $C$-universal, which contradicts the initial assumption. This now shows that $\hat{f}_\mathcal{K}(\lambda)$ must be strictly positive for all $\lambda \in \Lambda$. 

\subsubsection{Proof of Lemma \ref{lem:integralpd_four_comp}} \label{subsec:proof_c_lemma}

Since the kernel $\mathcal{K}$ is invariant,
$$
\mathcal{K}(g_1K,g_2K) = f_\mathcal{K}(g^{-1}_2g^{\phantom{-1}}_1\!\!\!\!K)
$$
for any $g_1, g_2 \in G$. Using (\ref{eq:compact_spherical}) and (\ref{eq:spherical_series}), 
$$
\mathcal{K}(g_1K,g_2K) = c_X\,\sum_{\lambda \in \Lambda} d_\lambda\hspace{0.02cm}\hat{f}_\mathcal{K}(\lambda)\hspace{0.03cm}\langle h_\lambda,U_\lambda(g^{-1}_2g^{\phantom{-1}}_1\!\!\!\!)\hspace{0.02cm}h_\lambda\rangle
$$
Since $U_\lambda$ is a unitary representation, this is
\begin{equation} \label{eq:compact_spectral_decompo}
\mathcal{K}(g_1K,g_2K) = c_X\,\sum_{\lambda \in \Lambda} d_\lambda\hspace{0.02cm}\hat{f}_\mathcal{K}(\lambda)\hspace{0.03cm}\langle U_\lambda(g_2)\hspace{0.02cm}h_\lambda,U_\lambda(g_1)\hspace{0.02cm}h_\lambda\rangle
\end{equation}
As explained after (\ref{eq:compact_fourier}), this can be written
\begin{equation} \label{eq:compact_spectral_decompo_bis}
\mathcal{K}(x,y) = c_X\,\sum_{\lambda \in \Lambda} d_\lambda\hspace{0.02cm}\hat{f}_\mathcal{K}(\lambda)\hspace{0.03cm}\langle (U_\lambda\hspace{0.02cm}h_\lambda)(y),(U_\lambda\hspace{0.02cm}h_\lambda)(x)\rangle    
\end{equation}
This is a uniformly convergent series, so the left-hand side of (\ref{eq:integralpd}) can be expressed as
$$
\int\int\,\mathcal{K}(x,y)\,\mu(dx)\mu(dy) = 
c_X\,\sum_{\lambda \in \Lambda} d_\lambda\hspace{0.02cm}\hat{f}_\mathcal{K}(\lambda)\left(\hspace{0.02cm}\int\int\langle (U_\lambda\hspace{0.02cm}h_\lambda)(y),(U_\lambda\hspace{0.02cm}h_\lambda)(x)\rangle\hspace{0.03cm}\mu(dx)\mu(dy)\right)
$$
Finally, by the definition of $\hat{\mu}(\lambda)$ in (\ref{eq:compact_fourier}),
\begin{align*}
    \text{LHS of (\ref{eq:integralpd})} & = 
c_X\,\sum_{\lambda \in \Lambda} d_\lambda\hspace{0.02cm}\hat{f}_\mathcal{K}(\lambda)\hspace{0.03cm} \langle\hat{\mu}(\lambda),\hat{\mu}(\lambda)\rangle \\
 &= 
c_X\,\sum_{\lambda \in \Lambda} d_\lambda\hspace{0.02cm}\hat{f}_\mathcal{K}(\lambda)\hspace{0.03cm} \Vert\hat{\mu}(\lambda)\Vert^2
\end{align*}
which is the required (\ref{eq:integralpd_four_comp}). 

\subsubsection{Proof of (\ref{eq:proof_comp_cov_1})} 

First, it will be shown that $\hat{\mu}(\lambda) = 0$ whenever $\lambda \neq \lambda_0\hspace{0.03cm}$. Indeed, by (\ref{eq:compact_fourier}) and the definition of $\mu$,
$$
\hat{\mu}(\lambda) = \int_X\,(U_\lambda\hspace{0.02cm}h_\lambda)(x)\hspace{0.03cm}m(x)\hspace{0.02cm}\mathrm{vol}(dx)
$$
This integral can be lifted to the group $G$. Specifically, 
if $dg$ denotes the Haar measure of $G$,  
\begin{equation} \label{eq:conv_comp_1}
\hat{\mu}(\lambda) = \left(c_X^{-1}\,\int_G\,U_\lambda(g)\hspace{0.02cm}m(gK)\hspace{0.02cm}dg\right)h_\lambda    
\end{equation}
where $c_X$ is the same as in (\ref{eq:spherical_series}). Now, recall that $m$ is given by a group convolution $m = \varphi^{\phantom{*}}_{\lambda_0} * \varphi^*_{\lambda_0}\hspace{0.03cm}$. Group convolutions obey the well-known convolution theorem~\cite{kirillov}
$$
\int_G\,U_\lambda(g)\hspace{0.02cm}m(gK)\hspace{0.02cm}dg = \left( \int_G\,U_\lambda(g)\hspace{0.02cm}\varphi^{\phantom{*}}_{\lambda_0}(gK)\hspace{0.02cm}dg\right)
\left( \int_G\,U_\lambda(g)\hspace{0.02cm}\varphi^*_{\lambda_0}(gK)\hspace{0.02cm}dg\right)
$$
From the definition (\ref{eq:compact_spherical}) of $\varphi_{\lambda_0}$ and the orthogonality relations for matrix coefficients of irreducible representations of $G$ (see~\cite{kirillov} or~\cite{sepanski}), the first integral on the right-hand side is zero for any $\lambda \neq \lambda_0\hspace{0.03cm}$. Since this is the case, the left-hand side is zero and $\hat{\mu}(\lambda) = 0$ by
(\ref{eq:conv_comp_1}). 

It remains to show that, as in (\ref{eq:proof_comp_cov_1}),
\begin{equation} \label{eq:sphere_coeff_m}
    \hat{\mu}(\lambda_0) = c(\lambda_0)h_{\lambda_0}    
\end{equation}
Using (\ref{eq:conv_comp_1}), it follows that for any $h \in V_{\lambda_0}$
\begin{equation} \label{eq:sphere_coeff_m_bis}
    \langle h,\hat{\mu}(\lambda_0)\rangle = c_X^{-1}\,\int_G\,\langle h,U_{\lambda_0}(g)\hspace{0.02cm}h_{\lambda_0}\rangle\hspace{0.03cm}m(gK)\hspace{0.02cm}dg
\end{equation}
Since the function $m$ is $K$-invariant, one has for any $k \in K$,
\begin{align}
\nonumber\langle h,\hat{\mu}(\lambda_0)\rangle & = c^{-1}_X\,\int_G\,\langle h,U_{\lambda_0}(k^{-1}g)\hspace{0.02cm}h_{\lambda_0}\rangle\hspace{0.03cm}m(gK)\hspace{0.02cm}dg \\ 
\label{eq:k_inv_comp} & = 
 c_X^{-1}\,\int_G\,\langle U_{\lambda_0}(k)\hspace{0.02cm}h,U_{\lambda_0}(g)\hspace{0.02cm}h_{\lambda_0}\rangle\hspace{0.03cm}m(gK)\hspace{0.02cm}dg    
\end{align}
where the second equality follows from the fact that $U_{\lambda_0}$ is a unitary representation.  

Note that (\ref{eq:k_inv_comp}) holds for any $k \in K$. Therefore, one may take its average over $k \in K$ (integrate both sides with respect to the Haar measure $dk$ of $K$). This yields,
$$
\langle h,\hat{\mu}(\lambda_0)\rangle = 
c_X^{-1}\,\int_G\int_K\,\langle U_{\lambda_0}(k)\hspace{0.02cm}h,U_{\lambda_0}(g)\hspace{0.02cm}h_{\lambda_0}\rangle\hspace{0.03cm}m(gK)\hspace{0.02cm}dk dg = \langle h^K,\hat{\mu}(\lambda_0)\rangle
$$
where the second equality is obtained by integrating first over $k$ and then over $g$. Specifically, $h^K$ is the average of $U_{\lambda_0}(k)\hspace{0.02cm}h$ over $k \in K$. This is just the orthogonal projection of $h$ onto the subspace of $K$-invariant vectors in $V_{\lambda_0}\hspace{0.03cm}$, a~one-dimensional subspace spanned by $h_{\lambda_0}$. It~therefore follows that $\hat{\mu}(\lambda_0) = c(\lambda_0)\hspace{0.02cm}h_{\lambda_0}$ for some $c(\lambda_0) \in \mathbb{C}$. Clearly,
$$
c(\lambda_0) = \langle\hat{\mu}(\lambda_0),h_{\lambda_0}\rangle
$$
But (\ref{eq:compact_spherical}) and (\ref{eq:sphere_coeff_m_bis}) then imply
$$
c^*(\lambda_0) = c^{-1}_X\,\int_G\,m(gK)\hspace{0.03cm}\varphi_{\lambda_0}(gK)\hspace{0.03cm}dg
$$
This is a real number because replacing $g$ by $g^{-1}$ leaves $m$ invariant and conjugates $\varphi_{\lambda_0}$ (recall~here that $m$ is real).

\section{Compact homogeneous spaces}\label{sec:compact_hom}
Theorem \ref{th:compact_univ} can be readily generalised to compact homogeneous spaces that are not necessarily symmetric. Here, $X = G/K$ where $G$ is a compact Lie group and $K$ a closed subgroup.

Consider the set $\Lambda$ of equivalence classes of irreducible representations $\lambda = (U_\lambda,V_\lambda)$ of $G$ which admit at least one non-zero $K$-invariant vector $h \in V_\lambda$ (that is $h \in V_\lambda$ with $U_\lambda(k)\hspace{0.02cm}h = h$ for all $k \in K$). For $\lambda \in \Lambda$, denote $W_\lambda$ the~subspace of $V_\lambda$ whose elements are all such $K$-invariant vectors.

Denote $d_\lambda$ and $r_\lambda$ the dimensions of $V_\lambda$ and $W_\lambda\hspace{0.03cm}$, respectively. In the previous section, $r_\lambda$ was equal to $1$ for all $\lambda \in \Lambda$, but this is no longer the case here. Choose some orthonormal basis $(h^\lambda_i\,;i=1,\ldots,r_\lambda)$ for each $W_\lambda$ and extend it to an~orthonormal basis $(h^\lambda_i\,;i=1,\ldots,d_\lambda)$ of $V_\lambda\hspace{0.03cm}$. Then, the matrix coefficients of the representation $U_\lambda$ will be denoted
\begin{equation}  \label{eq:hom_matrix}
    u_{ij}^\lambda(g) = \langle U_\lambda(g)\hspace{0.03cm}h^\lambda_j,h^\lambda_i\rangle \hspace{0.5cm} i,j = 1,\ldots,d_\lambda
\end{equation}
If $i,j \leq r_\lambda\hspace{0.03cm}$, $u_{ij}^\lambda(g)$ is constant on double cosets $KgK$, and can be identified with a $K$-invariant function $u_{ij}^\lambda(x)$ of $x \in X$. For continuous $K$-invariant functions $f:X \rightarrow \mathbb{C}$, the spherical series
 (\ref{eq:spherical_series}) here generalizes to (see~\cite{helgason2022groups}, Page 533)
\begin{equation}\label{eq:harmonic_series}
    f(x) = \sum_{\lambda \in \Lambda} d_\lambda \sum_{i=1}^{r_\lambda}\sum_{j=1}^{r_\lambda}
    \hspace{0.02cm}\hat{f}_{ij}(\lambda)\hspace{0.03cm}u_{ij}^\lambda(x) 
\end{equation}
where $\hat{f}_{ij}(\lambda)$ is given by ($dg$ denotes the Haar measure of $G$)
\begin{equation} \label{eq:harmonic_coeff}
    \hat{f}_{ij}(\lambda) = \int_G\,f(gK)\hspace{0.03cm}u^\lambda_{ji}(g^{-1})\hspace{0.03cm}dg
\end{equation}
In the following, let $\hat{f}(\lambda)$ denote the matrix with coefficients $\hat{f}_{ij}(\lambda)$ ($i,j = 1,\ldots,r_\lambda$). 

Recall Bochner's theorem for compact homogeneous spaces, which says that $\mathcal{K}:X \times X \rightarrow \mathbb{C}$ is continuous and positive-definite if and only if $\hat{f}_\mathcal{K}(\lambda)$ is a Hermitian positive-semidefinite matrix for each $\lambda \in \Lambda$ ($f_\mathcal{K}$ is the continuous $K$-invariant function $f_\mathcal{K}(x) = \mathcal{K}(x,K)$)~\cite{yaglom_second-order_1961}. It turns out that $\mathcal{K}$ is a universal kernel if and only if these matrices are (strictly) positive-definite. This is a direct generalization of Theorem \ref{th:compact_univ}.
\begin{theorem}\label{th:compact_hom_univ}
    A continuous, positive-definite invariant kernel $\mathcal{K}:X\times X \rightarrow \mathbb{R}$ is $C$-universal if and only if the matrix $\hat{f}_\mathcal{K}(\lambda)$ is Hermitian positive-definite for each $\lambda\in \Lambda$.
\end{theorem}
The proof of this theorem is mostly analogous to the proof of Theorem \ref{th:compact_univ}. The main difference is the proof of the only-if part, which will require using a complex version of Condition (\ref{eq:integralpd}). Specifically, a kernel $\mathcal{K}:X \times X \rightarrow \mathbb{R}$ satisfies (\ref{eq:integralpd}) for any non-zero finite signed measure $\mu$ on $X$, if and only if it satisfies
\begin{equation} \label{eq:integralpd_complex}
 \int_X\int_X\,\mathcal{K}(x,y)\mu(dx)\mu^*(dy) > 0
\end{equation}
for any non-zero finite complex measure $\mu$ on $X$ (where $^*$ denotes the complex conjugate)~\cite{sriperumbudur2011universality}\cite{guella2022}.

The proof requires a generalization of (\ref{eq:integralpd_four_comp}) from the previous section. This is the following (see the proof below)
\begin{equation}\label{eq:integralpd_four_comp_hom}
    \text{LHS of (\ref{eq:integralpd_complex})} = \sum_{\lambda\in\Lambda}d_\lambda\sum_{k=1}^{d_\lambda}
(\hat{\mu}_k(\lambda))^\dagger\hat{f}_{\mathcal{K}}(\lambda)\,(\hat{\mu}_k(\lambda))
\end{equation}
where $\hat{\mu}_k(\lambda) \in \mathbb{C}^{r_\lambda}$ is a column vector with components
\begin{equation} \label{eq:ulik}
    \hat{\mu}_{ki}(\lambda) = \int_X u_{ki}^\lambda(x)\hspace{0.03cm} \mu(dx) \hspace{0.5cm} i =1,\ldots, r_\lambda
\end{equation}
and the notation $u_{ki}^\lambda(x)$ is justified by the fact that $u_{ki}^\lambda(g)$ is constant on cosets $x = gK$ as long as $i = 1,\ldots, r_\lambda\hspace{0.03cm}$.

For the if part, similar to the proof of Theorem \ref{th:compact_univ}, it is enough to argue that any non-zero finite complex measure $\mu$ has at least one $\hat{\mu}_k(\lambda)$ different from zero. The corresponding term in (\ref{eq:integralpd_four_comp_hom}) has to be $> 0$, since $\hat{f}_\mathcal{K}(\lambda)$ is positive-definite.   

For the only-if part, suppose there exists $\lambda_0 \in \Lambda$ with $\hat{f}_\mathcal{K}(\lambda_0)$ not positive-definite. Pick some non-zero $v \in \mathbb{C}^{r_\lambda}$ such~that  $\hat{f}_\mathcal{K}(\lambda_0)\hspace{0.03cm}v = 0$. Let $\mu$ be a complex measure on $X$ with $\hat{\mu}_k(\lambda) = \delta_{\lambda\lambda_0}\hspace{0.02cm}v$ ($\delta_{\lambda\lambda_0} = 1$ if $\lambda = \lambda_0$ and $= 0$ otherwise). This is a non-zero complex measure which violates (\ref{eq:integralpd_complex}), as one may easily see by replacing into (\ref{eq:integralpd_four_comp_hom}). 

The existence of this $\mu$ is guaranteed by the fact that ($u^\lambda_{ki}\hspace{0.02cm};\,k=1,\ldots, d_\lambda\hspace{0.02cm}, i = 1, \ldots, r_\lambda$) are linearly independent and span a dense subset of $C(X,\mathbb{C})$~\cite{helgason2022groups} (Page 533), through a straightforward application of the Hahn-Banach theorem.

\subsubsection{Proof of (\ref{eq:integralpd_four_comp_hom})}
Recall that $f_\mathcal{K}(gK) = \mathcal{K}(gK,K)$ is continuous and constant on double cosets $KgK$. To simplify notation,\hfill\linebreak $f_\mathcal{K}$~is~here written as just $f$.
Then, by applying (\ref{eq:harmonic_series}) to this function, 
\begin{equation*}
    \begin{aligned}
   \mathcal{K}(g_1K,g_2K) = f(g^{-1}_2g^{\phantom{-1}}_1\!\!\!K) &=     \sum_{\lambda \in \Lambda} d_\lambda \sum_{i,j=1}^{r_\lambda}
    \hat{f}_{ij}(\lambda)\hspace{0.03cm}u_{ij}^\lambda(g^{-1}_2g^{\phantom{-1}}_1\!\!\!) \\
        &= \sum_{\lambda \in \Lambda} d_\lambda \sum_{k=1}^{d_\lambda}\sum_{i,j=1}^{r_\lambda}\hat{f}_{ij}(\lambda)\hspace{0.03cm} u_{ik}^\lambda(g_2^{-1})\hspace{0.03cm}u_{kj}^\lambda(g_1) \\
        &= \sum_{\lambda \in \Lambda} d_\lambda \sum_{k=1}^{d_\lambda}\sum_{i,j=1}^{r_\lambda}\hat{f}_{ij}(\lambda)\hspace{0.03cm} (u_{ki}^\lambda(g_2))^*\hspace{0.03cm}u_{kj}^\lambda(g_1) \\
    \end{aligned}
\end{equation*}
where the second and third equalities follow by using the fact that each $U_\lambda$ is a unitary representation. The series just obtained is uniformly convergent. Moreover, as explained after (\ref{eq:ulik}), matrix elements $u^\lambda_{ki}$ and $u^\lambda_{kj}$ (with $i,j = 1, \ldots r_\lambda$) are constant on the cosets $y = g_2K$ and $x = g_1K$. Replacing into the left-hand side of (\ref{eq:integralpd_complex}), it follows that
\begin{equation*}
    \begin{aligned}
            \int_X\int_X\,\mathcal{K}(x,y)\mu(dx)\mu^*(dy) 
            &= 
            \sum_{\lambda \in \Lambda} d_\lambda \sum_{k=1}^{d_\lambda}\sum_{i,j=1}^{r_\lambda}\hat{f}_{ij}(\lambda)\hspace{0.03cm} \int_X (u_{ki}^\lambda(y))^*\mu^*(dy)\int_X u_{kj}^\lambda(x)\mu(dx) \\
            &= \sum_{\lambda \in \Lambda} d_\lambda \sum_{k=1}^{d_\lambda}\sum_{i,j=1}^{r_\lambda}\hat{f}_{ij}(\lambda)\hspace{0.03cm}(\hat{\mu}_{ki}(\lambda))^*(\hat{\mu}_{kj}(\lambda))
    \end{aligned}
\end{equation*}
where the second equality is immediate from (\ref{eq:ulik}). Now, this last expression is the same as the right-hand side of (\ref{eq:integralpd_four_comp_hom}). 

\section{Non-compact symmetric spaces} \label{sec:noncompact}

\subsection{Statement of the main result} \label{subsec:noncompact}
In this section, $X = G/K$ will be a Riemannian symmetric space of the non-compact type. A~somewhat surprising result will be proved\,: if a $C_0$ invariant kernel is integrable (in a sense that will be made precise below), then it is $C_0$-universal. 

Here, a kernel is once again defined to be a Hermitian map $\mathcal{K}:X\times X \rightarrow \mathbb{C}$, which is said to be positive definite if for all $N\in\mathbb{N}$ and $x_1,\dots,x_N\in X$,
the Gram matrix $\big[\mathcal{K}(x_i,x_j)\big]_{i,j}$ is Hermitian positive semidefinite.
A continuous kernel $\mathcal{K}$ is called $C_0$ if $\mathcal{K}(x,y) \rightarrow 0$ when the Riemannian distance $d(x,y) \rightarrow \infty$ (with $x$ being fixed), and said to be invariant if $\mathcal{K}(g\hspace{0.02cm}x,g\hspace{0.02cm}y) = \mathcal{K}(x,y)$ for any $g \in G$ (as usual, $g\hspace{0.02cm} x$ denotes the action of $g \in G$ on $x \in G/K$ by multiplication from the left).

The kernel $\mathcal{K}$ is integrable  if the function $f_\mathcal{K}(x) = \mathcal{K}(x,K)$ belongs to $L^1(X,\mathrm{vol})$, where $\mathrm{vol}$ is the Riemannian volume measure of $X$. With these definitions, here is the main result.
\begin{theorem} \label{th:ncc0}
Let $\mathcal{K}:X \times X \rightarrow \mathbb{R}$ be a continuous, positive-definite and $C_0$ invariant kernel. If~$\mathcal{K}$ is integrable, then $\mathcal{K}$ is $C_0$-universal.
\end{theorem}
To prove this theorem, one needs to work with the spherical transform of the function $f_\mathcal{K}\hspace{0.03cm}$. Some basic facts about spherical functions and spherical transforms should be kept in mind. 

Note that $G$ is a non-compact real semisimple Lie group, so it has an Iwasawa decomposition $G = KAN$, where (of course) $K$ is compact, $A$ is abelian, and $N$ is nilpotent~\cite{helgason2022groups}. Denote $\mathfrak{a}$ the Lie algebra of $A$, and denote $\mathfrak{a}^*$ its dual (as a vector space). 

For each $\lambda \in \mathfrak{a}^*$, there is a spherical function $\varphi_\lambda:X \rightarrow \mathbb{C}$, given by the Harish-Chandra integral~\cite{helgason2022groups} ($dk$ denotes the normalised Haar measure of $K$), 
\begin{equation} \label{eq:sphericalfunction}
    \varphi_\lambda(gK) = \int_K\,e^{(\mathrm{i}\lambda - \rho)(a(gk))}\,dk
\end{equation}
with $\rho$ the half-sum of positive roots and $g = k\,\exp(a(g))\,n$ the Iwasawa decomposition of~$g \in G$. Then, each $\varphi_\lambda$ is~$K$-invariant (constant on double cosets $KgK$) and a bounded eigenfunction of the Laplace-Beltrami operator of $X$.  

If $f:X \rightarrow \mathbb{C}$ is $K$-invariant, and $f \in L^1(X,\mathrm{vol})$, its spherical transform is $\hat{f}:\mathfrak{a}^* \rightarrow \mathbb{C}$~\cite{helgason2022groups},
\begin{equation} \label{eq:FT}
    \hat{f}(\lambda) = \int_X\,f(x)\,\varphi_{-\lambda}(x)\hspace{0.02cm}\mathrm{vol}(dx) 
\end{equation}
and one has the inversion formula (assuming that the following integral converges absolutely),
\begin{equation} \label{eq:IFT}
    f(x) = c_X\,\int_{\mathfrak{a}^*}\,\hat{f}(\lambda)\,\varphi_\lambda(x)\hspace{0.02cm}|c(\lambda)|^{-2}d\lambda
\end{equation}
where $c_X$ is a normalising constant, and $c(\lambda)$ denotes the Harish-Chandra function~\cite{helgason2022groups}
(a~meromorphic function on $\mathfrak{a}^*_\mathbb{C}\hspace{0.03cm}$ -- the complexification of $\mathfrak{a}^*$), while $d\lambda$ is just the Lebesgue measure on $\mathfrak{a}^*$.  

The proof of Theorem \ref{th:ncc0} relies, first of all, on the $L^1$-Godement theorem~\cite{dacosta2024invariantkernelsriemanniansymmetric}. In particular, this theorem implies the following proposition.
\begin{proposition} \label{th:godement}
    Let $\mathcal{K}:X\times X \rightarrow \mathbb{C}$ be a continuous, positive-definite invariant kernel. If $\mathcal{K}$ is integrable, then the spherical transform $\hat{f}_\mathcal{K}$ is non-negative ($\hat{f}_\mathcal{K}(\lambda) \geq 0$ for all $\lambda \in \mathfrak{a}^*$) and integrable (with respect to the spectral measure $\sigma(d\lambda) = |c(\lambda)|^{-2}d\lambda$).
\end{proposition}
Note that if $\mathcal{K}$ is an invariant kernel, then $f_\mathcal{K}$ is a $K$-invariant function, so it makes sense to speak of its spherical transform. 

Another, crucial ingredient in the proof of Theorem \ref{th:ncc0} is the following claim about analyticity of spherical transforms~\cite{helgason2022groups} (Page 459). 

Let $C(\rho)$ denote the convex hull of the points $w\hspace{0.02cm}\rho \in \mathfrak{a}^*$ where $w$ runs through $W$, the Weyl group of the symmetric pair $(G,K)$, and $\rho$ is the half-sum of positive roots (as above, in (\ref{eq:sphericalfunction})). 
\begin{lemma} \label{lem:holomorphic}
    Let $f:X\rightarrow \mathbb{C}$ be $K$-invariant. If $f \in L^1(X,\mathrm{vol})$, then the spherical transform $\hat{f}:\mathfrak{a}^*\rightarrow \mathbb{C}$ extends to a holomorphic function on the tube domain $\mathfrak{a}^* + \mathrm{i}\hspace{0.02cm}C^o(\rho)$ in $\mathfrak{a}^*_\mathbb{C}$ (where $C^o(\rho)$ denotes the interior of $C(\rho)$).  
\end{lemma}

\subsection{Alternative form of (\ref{eq:integralpd})} \label{subsec:nctools}
Proposition \ref{th:godement} and Lemma \ref{lem:holomorphic} are not the whole story behind the proof of Theorem \ref{th:ncc0}. The final ingredient is an alternative expression of Condition (\ref{eq:integralpd}), which uses the concept of~group Fourier transform.

For each $\lambda \in \mathfrak{a}^*$, there is an irreducible unitary representation of $G$, say $U_\lambda\hspace{0.03cm}$. These representations all act in the same Hilbert space $L^2(K/M)$, where $M$ is the centraliser of $A$ in $K$, and $K/M$ is equipped with the measure induced by the Haar measure of $K$~\cite{anker}. Specifically, the action is the following
\begin{equation} \label{eq:sps}
    U_\lambda(g)\hspace{0.02cm}h(kM) = e^{-(\mathrm{i}\lambda + \rho)(a(g^{-1}k))}\,h(k(g^{-1}k)M)
\end{equation}
where $g^{-1}k$ has Iwasawa decomposition $k(g^{-1}k)\hspace{0.02cm}\exp(a(g^{-1}k))\hspace{0.02cm}n$.

In $L^2(K/M)$, there is a unique $K$-invariant unit vector, which is the constant function $h_0 = 1$. The spherical function $\varphi_\lambda$ in (\ref{eq:sphericalfunction}) is in fact~\cite{anker}
\begin{equation} \label{eq:sphericalrep}
      \varphi_\lambda(gK) = \langle h_0\hspace{0.02cm},U_\lambda(g)h_0\rangle
       \hspace{0.5cm}\text{(scalar product in $L^2(K/M)$)}
\end{equation}
For a finite signed measure $\mu$ on $G$, its group Fourier transform $\hat{\mu}$ associates to each $\lambda \in \mathfrak{a}^*$ a~bounded linear operator $\hat{\mu}(\lambda):L^2(K/M)\rightarrow L^2(K/M)$. This is the following~\cite{kirillov}, 
\begin{equation}
\label{eq:groupfour_measure}
     \hat{\mu}(\lambda) = \int_G\,U_\lambda(g)\hspace{0.02cm}\mu(dg)
\end{equation}
where the integral should be understood as a weak integral.

In Lemma \ref{lem:integralpd_four} below,
The group Fourier transform $\hat{\mu}$ will be used to express Condition (\ref{eq:integralpd}).
Note that measures $\mu$ which appear in (\ref{eq:integralpd}) are measures on $X = G/K$, not on the group $G$. However, a measure $\mu$ on $X$ can be lifted to a measure $\bar{\mu}$ on $G$, via the following construction
\begin{equation} \label{eq:averaging}
\int_G\,f(g)\hspace{0.02cm}\bar{\mu}(dg)  =   \int_X\,\bar{f}(x)\hspace{0.02cm}\mu(dx) \hspace{0.5cm}
\bar{f}(gK) = \int_K\,f(gk)\hspace{0.02cm}dk
\end{equation}
where $f:G \rightarrow \mathbb{R}$, while $\bar{f}:X\rightarrow \mathbb{R}$ is the ``right $K$-average" of $f$.  In the following, the measure $\mu$ on $X$ will not be distinguished from the measure $\bar{\mu}$ on $G$, and $\hat{\bar{\mu}}$ will just be denoted by $\hat{\mu}$.  
\begin{lemma}\label{lem:integralpd_four}
    Under the conditions of Proposition \ref{th:godement},
\begin{equation}
\label{eq:integralpd_four}   
\text{LHS of (\ref{eq:integralpd})} = c_X\,\int_{\mathfrak{a}^*}\,\hat{f}_\mathcal{K}(\lambda)\hspace{0.02cm}\Vert \hat{\mu}(\lambda)\hspace{0.02cm}h_0\Vert^2\hspace{0.03cm}\sigma(d\lambda)    
\end{equation}
for any finite signed measure $\mu$ on $X$. Here, the norm under the integral is that of $L^2(K/M)$, and $\sigma(d\lambda)$ denotes the spectral measure $\sigma(d\lambda) = |c(\lambda)|^{-2}d\lambda$.
\end{lemma}
The proof of Theorem \ref{th:ncc0} will require certain additional properties of $\hat{\mu}$. Above all,
\begin{proposition} \label{prop:injectivegroupfour}
    Let $\mu$ be a finite signed measure on $X$. The following properties hold
    \begin{itemize}
        \item[(a)] if $\hat{\mu}(\lambda)\hspace{0.02cm}h_0 = 0$ for all $\lambda \in \mathfrak{a}^*$, then $\mu$ is the zero measure.
        \item[(b)] $\Vert \hat{\mu}(\lambda)\hspace{0.02cm}h_0\Vert$ is a bounded continuous function of $\lambda \in \mathfrak{a}^*$.
    \end{itemize} 
\end{proposition}
Lemma \ref{lem:integralpd_four} and Proposition \ref{prop:injectivegroupfour} will be proved in the following section. Before going on, consider the connection between spherical transforms and group Fourier transforms. Whenever $f \in L^1(X,\mathrm{vol})$ is
 $K$-invariant as in (\ref{eq:FT}), $f$~can~be identified with a function $f:G\rightarrow \mathbb{C}$ which is constant on double cosets $KgK$. The group Fourier transform $\hat{f}$ of $f$ is just the group Fourier transform of the measure $\mu = f\hspace{0.02cm}dg$, where $dg$ denotes the Haar measure of $G$. 
 If~$(h_j\,;j=0,1,\ldots)$ is an orthonormal basis of $L^2(K/M)$ with $h_0 = 1$, consider the matrix coefficients
\begin{equation} \label{eq:groupfourij}
    \hat{f}_{ij}(\lambda) =
    \langle \hat{f}(\lambda)\hspace{0.02cm}h_j,h_i\rangle \hspace{0.5cm}  i,j = 0, 1, \ldots
\end{equation} 
\begin{proposition} \label{prop:grouptospherical}
In the above notation, the only non-zero matrix coefficient is $\hat{f}_{00}(\lambda)$, which is equal to the spherical transform in (\ref{eq:FT}).
\end{proposition}
This proposition is not needed below, so it will be proved separately in Section \ref{sec:grouptospherical}.

\subsection{Proof of the main result} \label{subsec:proof_ncc0}
Proposition \ref{th:godement} and Lemma \ref{lem:holomorphic} are known from~\cite{dacosta2024invariantkernelsriemanniansymmetric} and~\cite{helgason2022groups}. Admit, for now, that Lemma \ref{lem:integralpd_four} and Proposition \ref{prop:injectivegroupfour} are true.

To prove Theorem \ref{th:ncc0}, one has to show the left-hand side of (\ref{eq:integralpd}) is $> 0$ for any non-zero $\mu$. By Lemma \ref{lem:integralpd_four}, this is the same as the right-hand side of (\ref{eq:integralpd_four}). Assume, by way of contradiction,
\begin{equation} \label{eq:integralpd_contradiction}
\int_{\mathfrak{a}^*}\,\hat{f}_\mathcal{K}(\lambda)\hspace{0.02cm}\Vert \hat{\mu}(\lambda)\hspace{0.02cm}h_0\Vert^2\hspace{0.03cm}\sigma(d\lambda) = 0
\end{equation}
for some (as of now fixed) non-zero finite signed measure $\mu$. Recall first  $\sigma(d\lambda) = |c(\lambda)|^{-2}d\lambda$ where $c(\lambda)$ is the Harish-Chandra function. This function is given by the Gindikin-Karpelevich formula~\cite{helgason2022groups} (Page 447), which shows that $|c(\lambda)|^{-2} \neq 0$ as~long as $\lambda \in \mathfrak{a}^*$  and $\lambda \neq 0$ ($|c(\lambda)|^{-2}$ may still vanish for certain other $\lambda \in \mathfrak{a}^*_\mathbb{C}$). Accordingly, it follows from (\ref{eq:integralpd_contradiction}) that 
$\hat{f}_\mathcal{K}(\lambda)\Vert \hat{\mu}(\lambda)\hspace{0.02cm}h_0\Vert = 0$ for (Lebesgue) almost all $\lambda \in \mathfrak{a}^*$. However, since this is a continuous function (in~particular, due to Proposition \ref{prop:injectivegroupfour}), 
\begin{equation} \label{eq:integralpd_contradiction_bis}
    \hat{f}_\mathcal{K}(\lambda)\Vert \hat{\mu}(\lambda)\hspace{0.02cm}h_0\Vert = 0 \hspace{0.5cm} \text{for all $\lambda \in \mathfrak{a}^*$}
\end{equation}
Thus, $\Vert \hat{\mu}(\lambda)\hspace{0.02cm}h_0\Vert \neq 0$ only if $\lambda$ belongs to the vanishing set of $\hat{f}_\mathcal{K}(\lambda)$ in $\mathfrak{a}^*$ (the set of $\lambda \in \mathfrak{a}^*$ such that 
$\hat{f}_\mathcal{K}(\lambda) = 0$). By Lemma \ref{lem:holomorphic}, $\hat{f}_\mathcal{K}(\lambda)$ is a real-analytic function of $\lambda \in \mathfrak{a}^*$, its vanishing set therefore has empty interior (because an analytic function which vanishes on a non-empty open set vanishes identically~\cite{bochner}).

To conclude, recall that Proposition \ref{prop:injectivegroupfour} states $\Vert \hat{\mu}(\lambda)\hspace{0.02cm}h_0\Vert$ is continuous. A continuous function which vanishes on the complement of a set with empty interior must be identically zero. 
Accordingly, $\hat{\mu}(\lambda)\hspace{0.02cm}h_0 = 0$ identically on $\mathfrak{a}^*$, and~$\mu$~must be the zero measure, by Proposition \ref{prop:injectivegroupfour}.
This~shows that (\ref{eq:integralpd_contradiction}) cannot hold for any non-zero $\mu$. 

To complete the proof, it now only remains to prove Lemma \ref{lem:integralpd_four} and Proposition \ref{prop:injectivegroupfour}. 
\subsubsection{Proof of Lemma \ref{lem:integralpd_four}} \label{subsec:noncomp-lem_prooof} 
Proposition \ref{th:godement} entails $f_\mathcal{K}$ can be written down as an inverse spherical transform (as in (\ref{eq:IFT}))
\begin{equation} \label{eq:kernel_IFT1}
    f_\mathcal{K}(x) = 
c_X\,\int_{\mathfrak{a}^*}\,\hat{f}_\mathcal{K}(\lambda)\,\varphi_\lambda(x)\hspace{0.02cm}\sigma(d\lambda)
\end{equation}
By definition, $f_\mathcal{K}(x) = \mathcal{K}(x,K)$. Since the kernel $\mathcal{K}$ is invariant,
\begin{equation} \label{eq:kernel_IFT2}
    \mathcal{K}(g_1K,g_2K) = 
c_X\,\int_{\mathfrak{a}^*}\,\hat{f}_\mathcal{K}(\lambda)\,\varphi_\lambda(g^{-1}_2g^{\phantom{-1}}_1\!\!\!\!K)\hspace{0.02cm}\sigma(d\lambda)
\end{equation}
Then, by (\ref{eq:sphericalrep}),
\begin{align}
\nonumber 
\mathcal{K}(g_1K,g_2K)  & = 
c_X\,\int_{\mathfrak{a}^*}\,\hat{f}_\mathcal{K}(\lambda)\,\langle h_0,U_\lambda(g^{-1}_2g^{\phantom{-1}}_1\!\!\!\!)\hspace{0.02cm}h_0\rangle\hspace{0.02cm}\sigma(d\lambda) \\
\label{eq:kernel_IFT3} & = 
c_X\,\int_{\mathfrak{a}^*}\,\hat{f}_\mathcal{K}(\lambda)\,\langle U_\lambda(g_2)\hspace{0.02cm}h_0,U_\lambda(g_1)\hspace{0.02cm}h_0\rangle\hspace{0.02cm}\sigma(d\lambda)
\end{align}
where the second equality uses the fact that $U_\lambda$ is a unitary representation. To obtain (\ref{eq:integralpd_four}), it remains to use (\ref{eq:averaging}),
$$
    \int\int\,\mathcal{K}(x,y)\,\mu(dx)\mu(dy) = 
    \int\int\,\mathcal{K}(g_1K,g_2K)\,\bar{\mu}(dg_1)\bar{\mu}(dg_2)
$$
Replacing (\ref{eq:kernel_IFT3}) under the integral, and using Fubini's theorem
$$
\text{LHS of (\ref{eq:integralpd})} = 
c_X\,\int_{\mathfrak{a}^*}\hat{f}_\mathcal{K}(\lambda)\!\left(\int\int\,\langle U_\lambda(g_2)\hspace{0.02cm}h_0,U_\lambda(g_1)\hspace{0.02cm}h_0\rangle\hspace{0.02cm}\bar{\mu}(dg_1)\bar{\mu}(dg_2)\right)\sigma(d\lambda)
$$
By (\ref{eq:groupfour_measure}), the expression in parentheses is $\Vert \hat{\mu}(\lambda)\hspace{0.02cm}h_0\Vert^2$, just as in (\ref{eq:integralpd_four}). 

\subsubsection{Proof of Proposition \ref{prop:injectivegroupfour}} The second part (part (b)) is easier. To prove it, note that
$$
\Vert \hat{\mu}(\lambda)\hspace{0.02cm}h_0 \Vert^2 = \int\int\,\langle U_\lambda(g_2)\hspace{0.02cm}h_0,U_\lambda(g_1)\hspace{0.02cm}h_0\rangle\hspace{0.02cm}\bar{\mu}(dg_1)\bar{\mu}(dg_2) = 
\int\int\,\langle h_0,U_\lambda(g^{-1}_2g^{\phantom{-1}}_1\!\!\!\!)\hspace{0.02cm}h_0\rangle\hspace{0.02cm}\bar{\mu}(dg_1)\bar{\mu}(dg_2) 
$$
where the second equality uses the fact that $U_\lambda$ is a unitary representation (just as in the above proof of Lemma \ref{lem:integralpd_four}). Let $\nu$ be the image of the product measure $\bar{\mu}(dg_1)\bar{\mu}(dg_2)$ under the map $(g_1,g_2) \mapsto g^{-1}_2g^{\phantom{-1}}_1\!\!\!\!$ (this is again a finite signed measure). From the above,
$$
\Vert \hat{\mu}(\lambda)\hspace{0.02cm}h_0 \Vert^2 =
\int\,\langle h_0,U_\lambda(g)\hspace{0.02cm}h_0\rangle\hspace{0.02cm}\nu(dg) 
$$
Now, from (\ref{eq:sphericalrep}), it is clear that
$$
\Vert \hat{\mu}(\lambda)\hspace{0.02cm}h_0 \Vert^2 =
\int\, \varphi_\lambda(gK)\hspace{0.02cm}\nu(dg) 
$$
which is a continuous function of $\lambda \in \mathfrak{a}^*$, because $\varphi_\lambda(gK)$ is continuous and uniformly bounded ($|\varphi_\lambda(gK)| \leq 1$ for all $\lambda \in \mathfrak{a}^*$ and $g \in G$, by (\ref{eq:sphericalrep}) and the Cauchy-Schwarz inequality). 

Now, consider the proof of part (a). To begin, recall the heat kernel
$(\mathcal{P}_t\,;t > 0)$ of $X$ (see~\cite{gangolli})
\begin{equation} \label{eq:heat}
   \mathcal{P}_t(g_1K,g_2K) = c(t)\,\int_{\mathfrak{a}^*}\,\exp\!\left(-t\left((\lambda,\lambda) + (\rho,\rho)\right)\right)\varphi_\lambda(g^{-1}_2g^{\phantom{-1}}_1\!\!\!\!K)\hspace{0.03cm}\sigma(d\lambda)
\end{equation}
where $c(t)$ is a normalising constant and $(\cdot,\cdot)$ denotes the Killing form of the Lie algebra $\mathfrak{g}$ of $G$.\hfill\linebreak It is well-known that, for any continuous compactly supported $f:X\rightarrow \mathbb{R}$~\cite{davies},
\begin{equation} \label{eq:dirac}
f(x) = \lim_{t\rightarrow 0}\int_X\,f(y)\hspace{0.02cm}\mathcal{P}_t(x,y)\hspace{0.02cm}\mathrm{vol}(dy)    
\end{equation}
and the limit is uniform over $x$. Thus, for any finite-signed measure $\mu$ on $X$, it is easy to show
\begin{equation} \label{eq:diracmeasure}
    \int_X\,f(x)\hspace{0.02cm}\mu(dx) =
    \lim_{t\rightarrow 0}\int_X\,f(y)\left(\int_X\mathcal{P}_t(x,y)\mu(dx)\right)\mathrm{vol}(dy) 
\end{equation}
However, the heat kernel can be expressed as in (\ref{eq:kernel_IFT3}),
\begin{equation} \label{eq:heat_2}
   \mathcal{P}_t(g_1K,g_2K) = c(t)\,\int_{\mathfrak{a}^*}\,
\hat{f}_t(\lambda)\,
\langle U_\lambda(g_2)\hspace{0.02cm}h_0,
U_\lambda(g_1)\hspace{0.02cm}h_0\rangle   
\hspace{0.03cm}\sigma(d\lambda)    
\end{equation}
where $\hat{f}_t(\lambda) =    \exp\!\left(-t\left((\lambda,\lambda) + (\rho,\rho)\right)\right)$ as in (\ref{eq:heat}). Then, from (\ref{eq:averaging}),
$$
\int_X\,\mathcal{P}_t(x,y)\hspace{0.02cm}\mu(dx) = \int_G\,\mathcal{P}_t(g_1K,g_2K)\hspace{0.02cm}\,\bar{\mu}(dg_1)
$$
where $y = g_2K$. Moreover, from (\ref{eq:heat_2}), using Fubini's theorem,
$$
\int_X\,\mathcal{P}_t(x,y)\hspace{0.02cm}\mu(dx) = c(t)\,\int_{\mathfrak{a}^*}\,
\hat{f}_t(\lambda)\,
\langle U_\lambda(g_2)\hspace{0.02cm}h_0,
\hat{\mu}(\lambda)\hspace{0.02cm}h_0\rangle   
\hspace{0.03cm}\sigma(d\lambda) 
$$
Accordingly, if $\hat{\mu}(\lambda)\hspace{0.02cm}h_0$ for all $\lambda \in \mathfrak{a}^*$,
$$
\int_X\,\mathcal{P}_t(x,y)\hspace{0.02cm}\mu(dx) = 0
$$
for any $y \in X$, and (\ref{eq:diracmeasure}) implies that 
$$
\int_X\,f(x)\hspace{0.02cm}\mu(dx) = 0
$$
for any continuous compactly supported $f:X \rightarrow \mathbb{R}$. The same is then true for any bounded continuous $f$, which means that $\mu$ is the zero measure. 

\section{Revisiting the group Fourier transform} \label{sec:grouptospherical}
The group Fourier transform (\ref{eq:groupfour_measure}) is at the heart of harmonic analysis on any kind of group. Here, the aim is to explain that it takes on a specific, concrete form when applied to symmetric spaces of the non-compact type. 
In~particular, there is a  convolution theorem, inversion formula, and Plancherel formula, which are quite specific to this symmetric space setting.

To begin, consider the following proposition, which includes Proposition \ref{prop:grouptospherical} as a special case. A finite signed measure $\mu$ on $G$ is called $K$-invariant on the right (respectively, on the left), if~it is equal to its own image under any right translation $g \mapsto gk$ where $k \in K$ (respectively, any left translation $g \mapsto kg$). 
\begin{proposition} \label{prop:grouptospherical_general}
In the notation of (\ref{eq:groupfour_measure}), consider the matrix coefficients
\begin{equation} \label{eq:groupfour_ij_general}
   \hat{\mu}_{ij}(\lambda) =
    \langle \hat{\mu}(\lambda)\hspace{0.02cm}h_j,h_i\rangle  \hspace{0.5cm} i,j = 0,1, \ldots     
\end{equation}
If $\mu$ is $K$-invariant on the right (respectively, on the left), then $\hat{\mu}_{ij}(\lambda) = 0$ whenever $j > 0$ (respectively, whenever $i > 0$).
\end{proposition}
This proposition will be proved below. It has Proposition \ref{prop:grouptospherical} as an immediate consequence. To see this, note that $f \in L^1(X,\mathrm{vol})$  can be identified with $\tilde{f} \in L^1(G,dg)$ where $\tilde{f}(g) = f(gK)$. The group Fourier transform $\hat{f}(\lambda)$ is just the group Fourier transform of the measure $\mu = \tilde{f}dg$. From the very definition of $\tilde{f}$, this measure is $K$-invariant on the right. Proposition \ref{prop:grouptospherical_general} implies $\hat{f}_{ij}(\lambda) = 0$ if $j > 0$. If, in addition, $f$ is $K$-invariant, then $\mu$ is also $K$-invariant on the left, and (by a second application of Proposition \ref{prop:grouptospherical_general}) $\hat{f}_{ij}(\lambda) = 0$ if $i+j > 0$. This only leaves out $\hat{f}_{00}(\lambda)$.\hfill\linebreak By (\ref{eq:groupfour_measure}) and (\ref{eq:groupfour_ij_general}), 
$$
\hat{f}_{00}(\lambda) = \int_G\,\tilde{f}(g)\,\langle U_\lambda(g)\hspace{0.02cm}h_0,h_0\rangle\hspace{0.03cm}dg
$$
Recalling (\ref{eq:sphericalrep}), and using the fact that $dg$ is related to $\mathrm{vol}$ in the same way as $\bar{\mu}$ to $\mu$ in (\ref{eq:averaging})
\cite{helgason2022groups} (Page 91),
$$
\hat{f}_{00}(\lambda) = \int_G\,\tilde{f}(g)\,\varphi^*_\lambda(gK)\hspace{0.03cm}dg = \int_X\,f(x)\,\varphi^*_\lambda(x)\hspace{0.03cm}\mathrm{vol}(dx)
$$
where $^*$ denotes the complex conjugate. This is the same as the spherical transform (\ref{eq:FT}), because 
$\varphi^*_\lambda(x) = \varphi_{-\lambda}(x)$ for $\lambda \in \mathfrak{a}^*$ (see  (\ref{eq:sphericalfunction})). This shows that Proposition \ref{prop:grouptospherical} is indeed true. 

From the above discussion of Proposition \ref{prop:grouptospherical}, it appears that the group Fourier transform $\hat{f}(\lambda)$ of $f \in L^1(X,\mathrm{vol})$ is completely given by the coefficients $\hat{f}_i(\lambda) = \hat{f}_{i0}(\lambda)$ (where $i = 0,1, \ldots)$. In fact, these can be used to reconstruct the function $f$.


\begin{proposition} \label{prop:fourier_inversion}
Consider the functions $\varphi^i_\lambda(gK) = \langle h_i,U_\lambda(g)\hspace{0.02cm}h_0\rangle$ ($\lambda \in \mathfrak{a}^*$ and $i = 0,1,\ldots$). If $f \in L^1(X,\mathrm{vol})$ is continuous, and if
\begin{equation} \label{eq:inversion_summability}
\sum^\infty_{i=0}\,\int_{\mathfrak{a}^*}|\hat{f}_i(\lambda)|\hspace{0.02cm}\sigma(d\lambda) < \infty
\end{equation}
then, for all $x \in X$,
\begin{equation} \label{eq:groupfour_inversion}
    f(x) = c_X\,\sum^\infty_{i=0}\,\int_{\mathfrak{a}^*}\,\hat{f}_i(\lambda)\hspace{0.02cm}\varphi^i_\lambda(x)\hspace{0.03cm}\sigma(d\lambda)
\end{equation}
where $c_X$ is a normalising constant (the same as in (\ref{eq:IFT})). 
\end{proposition}
Note that the inversion formula (\ref{eq:groupfour_inversion}) reduces to the spherical transform inversion formula (\ref{eq:IFT}) when $f$ is $K$-invariant. This follows by a straightforward application of Proposition \ref{prop:grouptospherical}. 

Consider now the convolution theorem. For $t \in L^1(G,dg)$ and $f \in L^1(X,\mathrm{vol})$, define
\begin{equation} \label{eq:convdef_nc}
    (t*f)(x) = \int_G\,t(h)\hspace{0.02cm}f(h^{-1}x)\hspace{0.03cm}dh
\end{equation}
This is again a function in $L^1(X,\mathrm{vol})$, the convolution of $t$ and $f$.
\begin{proposition} \label{prop:convth_nc}
    For $t \in L^1(G,dg)$ and $f \in L^1(X,\mathrm{vol})$, define $t*f \in L^1(X,\mathrm{vol})$ as in (\ref{eq:convdef_nc}). Then,
    \begin{equation} \label{eq:convth_nc}
        (\widehat{t*f})_i(\lambda) = \sum^\infty_{j=0}\,\hat{t}_{ij}(\lambda)\hspace{0.02cm}
        \hat{f}_{j}(\lambda)
    \end{equation}
\end{proposition}
Using this proposition, one may recover a Plancherel theorem.
\begin{proposition} \label{prop:planc_nc}
    For $f \in L^1(X,\mathrm{vol}) \cap L^2(X,\mathrm{vol})$,
\begin{equation} \label{eq:plancherel_nc}
    \int_X\,|f(x)|^2\hspace{0.03cm}\mathrm{vol}(dx) = 
c_X\sum^\infty_{i=0}\,\int_{\mathfrak{a}^*}\,|\hat{f}_i(\lambda)|^2\hspace{0.03cm}\sigma(d\lambda)
\end{equation}
\end{proposition}
Propositions \ref{prop:fourier_inversion}--\ref{prop:planc_nc}
show that the group Fourier transform provides a suitable concept of ``Fourier transform" 
for~functions on symmetric spaces of the non-compact type. However, the established ``Fourier transform" for these functions is the Fourier-Helgason transform of~\cite{helgasonthirdbook}. How are these two related, and how are they different?

For $f \in L^1(X,\mathrm{vol})$, consider the integral
\begin{equation} \label{eq:prehelgason}
    F(\lambda,kM) = \int_G\,\tilde{f}(g)\left(U_\lambda(g)\hspace{0.02cm}h_0\right)\!(kM)\hspace{0.03cm}dg
\end{equation}
where $\tilde{f}(g) = f(gK)$, $U_\lambda$ is given by (\ref{eq:sps}) and $h_0 \in L^2(K/M)$ is the constant function $h_0 = 1$. If $F(\lambda,kM)$ is~understood as a weak integral, then it is defined as follows
\begin{align}
\nonumber F(\lambda,kM) &= 
    \sum^\infty_{i=0}\left(\hspace{0.02cm}\int_G\,\tilde{f}(g)\hspace{0.03cm}\langle U_\lambda(g)\hspace{0.02cm}h_0,h_i\rangle\hspace{0.03cm}dg\right)h_i(kM)  \\
    \label{eq:weak_helgason}
    &=  \sum^\infty_{i=0}\,\hat{f}_i(\lambda)\hspace{0.02cm}h_i(kM)
\end{align} 
where the coefficients $\hat{f}_i(\lambda)$ determine the group Fourier transform of $f$, as in Propositions \ref{prop:fourier_inversion}--\ref{prop:planc_nc}. On the other hand, if $F(\lambda,kM)$ is understood as a strong integral, then (\ref{eq:sps}) implies  
\begin{equation} \label{eq:strong_helgason}
    F(\lambda,kM) = \int_G\,\tilde{f}(g)\hspace{0.02cm}
    e^{-(\mathrm{i}\lambda + \rho)(a(g^{-1}k))}\hspace{0.02cm}dg
\end{equation}
which is essentially the Fourier-Helgason transform~\cite{helgasonthirdbook} (see Page 223). In the present work, the~reason for using the group Fourier transform instead of the Fourier-Helgason transform is that $F(\lambda,kM)$ exists as a weak integral as soon as $f \in L^1(X,\mathrm{vol})$, whereas the strong integral (\ref{eq:strong_helgason}) may diverge for $f \in L^1(X,\mathrm{vol})$. In the proof of Theorem \ref{th:ncc0}, the group Fourier transform had to be applied to finite signed measures $\mu$ (see Lemma \ref{lem:integralpd_four} and Proposition \ref{prop:injectivegroupfour}), which~include the special case $\mu = \tilde{f}dg$ with $f \in L^1(X,\mathrm{vol})$. 

Before moving on to the proofs of Propositions \ref{prop:grouptospherical_general}--\ref{prop:planc_nc}, a further remark about group Fourier transforms is in order. The representations of the form $U_\lambda$ in (\ref{eq:sps}) are only a (proper) subset of all the irreducible unitary representations of $G$ --- this subset is known as the spherical principal series~\cite{knappreps}. Of course, the group Fourier transform can be defined for any irreducible unitary representation $U$
 (not necessarily of the form (\ref{eq:sps})), just by putting $U$ instead of $U_\lambda$ in (\ref{eq:groupfour_measure}). In the present context, this is not necessary, since spherical principal series representations are already sufficient for Fourier analysis of integrable functions on the symmetric space $X = G/K$. 

Finally, here are the proofs of Propositions \ref{prop:grouptospherical_general} through \ref{prop:planc_nc}. They are mostly rather straightforward and use standard ideas of harmonic analysis. Therefore, they will be presented in a more succinct form, merely pointing out the main steps of each proof. 

\subsection{Proofs of Propositions \ref{prop:grouptospherical_general}--\ref{prop:planc_nc}} \label{subsec:groupfour_proofs}

\subsubsection{Proof of Proposition \ref{prop:grouptospherical_general}} Note that the proofs of the statements concerning $K$-invariance on the right and on the left are essentially identical. Consider the former. From (\ref{eq:groupfour_measure}) and (\ref{eq:groupfour_ij_general}),
$$
\hat{\mu}_{ij}(\lambda) = \int_G\,\langle U_\lambda(g)\hspace{0.02cm}h_j,h_i\rangle\hspace{0.02cm}\mu(dg)
$$
If $\mu$ is $K$-invariant on the right, this equality continues to hold when $U_\lambda(g)$ is replaced with $U_\lambda(gk)$, for any $k \in K$, under the integral. Integrating with respect to the normalised Haar measure $dk$, it follows that 
$$
\hat{\mu}_{ij}(\lambda) = \int_G\left(\hspace{0.02cm}\int_K\,\langle U_\lambda(gk)\hspace{0.02cm}h_j,h_i\rangle\hspace{0.02cm}dk\right)\mu(dg)
$$
Using the fact that $U_\lambda$ is a representation (a group homomorphism), the integral with respect to $dk$ becomes
$$
\langle U_\lambda(g)\hspace{0.02cm}h^K_j,h_i\rangle \;\;\text{ where }\;\; h^K_j = \int_K\,U_\lambda(k)\hspace{0.02cm}h_j\hspace{0.03cm}dk
$$
The linear map which takes $h_j$ to $h^K_j$ is the orthogonal projection (in $L^2(K/M)$) onto the subspace of $K$-invariant functions. This is a one-dimensional subspace, generated by the constant function $h_0\hspace{0.03cm}$. In other words, $h^K_j = h_0$ if~$j = 0$~and $=0$ if $j > 0$, and this immediately yields the required statement for $K$-invariance on the right. 

\subsubsection{Proof of Proposition \ref{prop:fourier_inversion}} To begin, recall the heat kernel
$(\mathcal{P}_t\,;t > 0)$, given by (\ref{eq:heat}). If $f \in L^1(X,\mathrm{vol})$ is continuous, then~\cite{davies}
$$
  f(x) = \lim_{t\rightarrow 0}\,\int_X\,f(y)\hspace{0.02cm}\mathcal{P}_t(x,y)\hspace{0.02cm}\mathrm{vol}(dy)    
$$
where the limit holds for each $x \in X$. Then, using (\ref{eq:heat_2}) and changing the order of integration
\begin{equation} \label{eq:l1_heat_spectral}
f(g_1K) = \lim_{t\rightarrow 0}\,
c(t)\,\int_{\mathfrak{a}^*}\,
\hat{f}_t(\lambda)\left(\hspace{0.02cm}\int_G\,\tilde{f}(g_2)\hspace{0.02cm}
\langle U_\lambda(g_2)\hspace{0.02cm}h_0,
U_\lambda(g_1)\hspace{0.02cm}h_0\rangle\hspace{0.03cm}dg_2\right)\sigma(d\lambda)  
\end{equation}
Here, the scalar product can be written,
$$
\langle U_\lambda(g_2)\hspace{0.02cm}h_0,
U_\lambda(g_1)\hspace{0.02cm}h_0\rangle = 
\sum^\infty_{i=0}\, \langle U_\lambda(g_2)\hspace{0.02cm}h_0,h_i\rangle\hspace{0.03cm}\langle h_i,U_\lambda(g_1)\hspace{0.03cm}h_0\rangle
$$
so that the interior integral, with respect to $dg_2$, can be expressed
$$
\sum^\infty_{i=0}\left(\int_G\,\tilde{f}(g_2)\langle U_\lambda(g_2)\hspace{0.02cm}h_0,h_i\rangle\hspace{0.03cm}dg_2\right)\langle h_i,U_\lambda(g_1)\hspace{0.02cm}h_0\rangle = \sum^\infty_{i=0}\,\hat{f}_i(\lambda)\hspace{0.03cm}\varphi^i_\lambda(g_1K) 
$$
Replacing this into (\ref{eq:l1_heat_spectral}), it follows that
\begin{equation} \label{eq:sum_inversion}
    f(gK) = \lim_{t\rightarrow 0}\,
c(t)\,\int_{\mathfrak{a}^*}\,
\hat{f}_t(\lambda)\left( \sum^\infty_{i=0}\,\hat{f}_i(\lambda)\hspace{0.03cm}\varphi^i_\lambda(gK) \right)\sigma(d\lambda)  
\end{equation}
To conclude the proof of (\ref{eq:groupfour_inversion}), it is now enough to apply the dominated convergence theorem. This is possible due to Condition 
(\ref{eq:inversion_summability}), and the fact that $|\varphi^i_\lambda(gK)| \leq 1$ by the definition of $\varphi^i_\lambda$ (using the Cauchy-Schwarz inequality). In~this~way, (\ref{eq:groupfour_inversion}) follows by taking the limit under the integral (\ref{eq:sum_inversion}) and noting that $c(t)\hspace{0.02cm}\hat{f}_t(\lambda)$ converges to the positive constant $c_X$ as $t \rightarrow 0$. Then,
$$
    f(gK) = 
c_X\,\int_{\mathfrak{a}^*}\,\sum^\infty_{i=0}\,\hat{f}_i(\lambda)\hspace{0.03cm}\varphi^i_\lambda(gK)\sigma(d\lambda)  
$$
where the integral and the sum can be interchanged, again due to (\ref{eq:inversion_summability}). 

\subsubsection{Proof of Proposition \ref{prop:convth_nc}} According to (\ref{eq:convdef_nc}),
\begin{align*}
(\widehat{t*f})_i(\lambda) & = \int_G\left(\hspace{0.02cm}\int_G\,t(h)\hspace{0.02cm}\tilde{f}(h^{-1}g)\hspace{0.02cm}dh\right)\langle U_\lambda(g)\hspace{0.02cm}h_0,h_i\rangle\hspace{0.03cm}dg    \\[0.15cm]
& = \int_G\int_G\,t(h)\hspace{0.02cm}\tilde{f}(h^{-1}g)\hspace{0.02cm}\langle U_\lambda(g)\hspace{0.02cm}h_0,h_i\rangle\hspace{0.03cm}dh\hspace{0.02cm}dg
\end{align*}
Introducing the new variables $g_1 = h^{-1}g$ and $g_2 = h$, 
$$
(\widehat{t*f})_i(\lambda) = \int_G\int_G\,t(g_2)\hspace{0.02cm}\tilde{f}(g_1)\hspace{0.02cm}\langle U_\lambda(g_2\hspace{0.02cm}g_1)\hspace{0.02cm}h_0,h_i\rangle\hspace{0.03cm}dg_1\hspace{0.02cm}dg_2
$$
Since $U_\lambda$ is a unitary representation, this is the same as
$$
(\widehat{t*f})_i(\lambda) = \int_G\int_G\,t(g_2)\hspace{0.02cm}\tilde{f}(g_1)\hspace{0.02cm}\langle U_\lambda(g^{\phantom{-1}}_1\!\!\!\!)\hspace{0.02cm}h_0,U_\lambda(g^{-1}_2)\hspace{0.02cm}h_i\rangle\hspace{0.03cm}dg_1\hspace{0.02cm}dg_2
$$
However, by expanding the scalar product in the orthonormal basis $(h_j)$,
\begin{align*}
(\widehat{t*f})_i(\lambda) & =
\sum^\infty_{j=0}\,\int_G\int_G\,t(g_2)\hspace{0.02cm}\tilde{f}(g_1)\hspace{0.03cm}\langle U_\lambda(g^{\phantom{-1}}_1\!\!\!\!)\hspace{0.02cm}h_0,h_j\rangle\hspace{0.03cm}\langle h_j,U_\lambda(g^{-1}_2)\hspace{0.02cm}h_i\rangle\hspace{0.03cm}dg_1\hspace{0.02cm}dg_2  \\[0.15cm]
& =
\sum^\infty_{j=0}\left(\hspace{0.02cm}\int_G\,t(g_2)\hspace{0.02cm}
\langle h_j,U_\lambda(g^{-1}_2)\hspace{0.02cm}h_i\rangle\hspace{0.03cm}dg_2\right)\left(\hspace{0.02cm}\int_G\,\tilde{f}(g_1)\hspace{0.03cm}\langle U_\lambda(g^{\phantom{-1}}_1\!\!\!\!)\hspace{0.02cm}h_0,h_j\rangle\hspace{0.03cm}dg_1\right)
\end{align*}
Using one last time the fact that $U_\lambda$ is unitary,
$$
(\widehat{t*f})_i(\lambda)  =
\sum^\infty_{j=0}\left(\hspace{0.02cm}\int_G\,t(g_2)\hspace{0.02cm}
\langle U_\lambda(g_2)\hspace{0.02cm}h_j,h_i\rangle\hspace{0.03cm}dg_2\right)\left(\hspace{0.02cm}\int_G\,\tilde{f}(g_1)\hspace{0.03cm}\langle U_\lambda(g_1)\hspace{0.02cm}h_0,h_j\rangle\hspace{0.03cm}dg_1\right)
$$
From the definition of $\hat{t}_{ij}(\lambda)$ and $\hat{f}_j(\lambda)$, this is exactly (\ref{eq:convth_nc}).

\subsubsection{Proof of Proposition \ref{prop:planc_nc}} The main steps of the proof are indicated without additional detail. The idea is to consider the convolution (\ref{eq:convdef_nc}), with $t(g) = f^\dagger(g)$ where $f^\dagger(g) = \tilde{f}^*(g^{-1})$ ($^*$ denotes the complex conjugate). By a straightforward calculation
\begin{equation} \label{eq:planc_convolution}
    (f^\dagger*f)(K) = \int_X\,|f(x)|^2\hspace{0.03cm}\mathrm{vol}(dx)
\end{equation}
where $(f^\dagger*f)(K)$ is the value of $f^\dagger*f$ at the coset $K \in G/K$. This is the left-hand side of (\ref{eq:plancherel_nc}). The proof of (\ref{eq:plancherel_nc}) relies on applying the spherical inversion formula (\ref{eq:IFT}) to $f^\dagger*f$. Indeed, it is possible to show that $f^\dagger*f:X\rightarrow \mathbb{C}$ is $K$-invariant and also $f^\dagger*f \in L^1(X,\mathrm{vol})$. Moreover, (\ref{eq:IFT}) holds true for this function, and reads
$$
(f^\dagger*f)(K) = c_X\,\int_{\mathfrak{a}^*}\,
\widehat{(f^\dagger*f)}(\lambda)\hspace{0.02cm}\varphi_\lambda(K)\hspace{0.03cm}\sigma(d\lambda)= c_X\,\int_{\mathfrak{a}^*}\,
\widehat{(f^\dagger*f)}(\lambda)\hspace{0.03cm}\sigma(d\lambda)
$$
where the second equality follows because $\varphi_\lambda(K) = 1$ for all $\lambda \in \mathfrak{a}^*$. It only remains to prove
\begin{equation} \label{eq:plancherelproof}
    \widehat{(f^\dagger*f)}(\lambda) = \sum^\infty_{i=0}\,|\hat{f}_i(\lambda)|^2
\end{equation}
However, this follows from Proposition \ref{prop:convth_nc}, after showing that $\widehat{f^\dagger}_{ij}(\lambda) = \hat{f}^*_j(\lambda)$ 
if $i = 0$ and $\widehat{f^\dagger}_{ij}(\lambda) =0$ if $i > 0$. By~Proposition \ref{prop:grouptospherical_general}, it is immediate that $\widehat{f^\dagger}_{ij}(\lambda) = 0$ for $i > 0$. On the other hand, 
$$   
\widehat{f^\dagger}_{0j}(\lambda) = \int_G\,\tilde{f}^*(g^{-1})\hspace{0.02cm}\langle U_\lambda(g)\hspace{0.02cm}h_j,h_0\rangle   \hspace{0.03cm}dg = 
\int_G\,\tilde{f}^*(g)\hspace{0.02cm}\langle U_\lambda(g^{-1})\hspace{0.02cm}h_j,h_0   \rangle\hspace{0.03cm}dg
$$
where the second equality follows after the change of variables $g \mapsto g^{-1}$. Finally, using the fact that $U_\lambda$ is unitary, 
it~is~easy to show that this is exactly the same as $\hat{f}^*_j(\lambda)$.      

\section{Examples and applications} \label{sec:examples} 

\subsection{Compact symmetric spaces}
Theorem \ref{th:compact_univ} says that continuous, positive-definite universal kernels $\mathcal{K}$ are precisely those that have no coefficients $\hat{f}_\mathcal{K}(\lambda)$ equal to zero in the spherical series (\ref{eq:spherical_series}) of $f_\mathcal{K}\hspace{0.03cm}$. The special cases of spheres and compact Lie groups are now considered. 

\subsubsection{Spheres} Let $M = S^n$ be the unit sphere in $\mathbb{R}^{n+1}$. The invariance group is the special orthogonal group $G = SO(n+1)$, and~any invariant kernel is of the form $\mathcal{K}(x,y) = g(x\cdot y)$ where $g$ is a continuous function and $x\cdot y$ denotes the Euclidean scalar product. 

The set $\Lambda$ of spherical irreducible representations is identified with the set of non-negative integers $\lambda$. To each $\lambda$ corresponds the regular representation of $SO(n+1)$ in the space of homogeneous harmonic polynomials of degree $\lambda$ defined on $\mathbb{R}^{n+1}$~\cite{vilenkin}. Moreover, the spherical function $\varphi_\lambda$ in (\ref{eq:compact_spherical}) is given by
\begin{equation} \label{eq:gegenbauer}
    \varphi_\lambda(x) = C^{(n-1)/2}_\lambda(x_{n+1})/C^{(n-1)/2}_\lambda(1)
\end{equation}
where $C^\alpha_l$ denotes the Gegenbauer polynomial of order $\alpha$ and degree $l$, and $x = (x_1,\ldots,x_{n+1})$. The spherical series (\ref{eq:spherical_series}) is therefore just an expansion in terms of Gegenbauer polynomials. For any invariant kernel $\mathcal{K}$, this reads
\begin{align} 
\nonumber    \mathcal{K}(x,y) = & \\ 
\label{eq:spherseries}
g(x\cdot y) = &\sum^\infty_{\lambda = 0} \hat{g}(\lambda)\hspace{0.03cm}C^{(n-1)/2}_\lambda(x\cdot y)
\end{align}
and $\mathcal{K}$ is universal exactly when all the coefficients $\hat{g}(\lambda)$ are strictly positive. Thus, $\mathcal{K}$ cannot be universal if $g$ is a polynomial, or if $g$ is an even function (as $\hat{g}(\lambda)$ would vanish for odd $\lambda$). On the other hand, the heat kernel on the sphere $S^n$ has the following expansion~\cite{vilenkin}
\begin{equation} \label{eq:sphere_heat}
    \mathcal{K}_t(x,y) = 
    \sum^\infty_{\lambda = 0} ((2\lambda +n-1)/(n-1))\hspace{0.02cm}e^{-t\hspace{0.02cm}\lambda(\lambda+n-1)}\hspace{0.03cm}C^{(n-1)/2}_\lambda(x\cdot y)
\end{equation}
where $t > 0$ is the time parameter. It is therefore clearly universal for any $t > 0$.

Another recipe for a universal kernel follows from the generating function of Gegenbauer polynomials~\cite{vilenkin},
$$
(1-2ts+t^2)^{-\alpha} = \sum^\infty_{l=0}\,C^\alpha_l(s)\hspace{0.02cm}t^l
\hspace{0.5cm} -1 \leq t \leq 1
$$
Based on this expansion, consider the kernel
\begin{equation} \label{eq:newtonianpotential}
    \mathcal{K}(x,y) = (1-2t\hspace{0.03cm}x\cdot y + t^2)^{-(n-1)/2} \hspace{0.5cm} 0 < t \leq 1
\end{equation}
This has the coefficients $\hat{g}(\lambda) = t^\lambda > 0$ and is therefore universal.

\subsubsection{Compact Lie groups} Let $M = G$ be a compact semisimple Lie group, viewed as a quotient $(G\times G)/\Delta$ where $\Delta$ is the diagonal subgroup of $G \times G$. Any invariant kernel is then of the form $\mathcal{K}(x,y) = f_\mathcal{K}(y^{-1}x)$ where $f_\mathcal{K}:G \rightarrow \mathbb{C}$ is a class function, $f_\mathcal{K}(g\hspace{0.02cm}x\hspace{0.02cm}g^{-1}) = f_\mathcal{K}(x)$ for $g,x \in G$. Moreover, the spherical series (\ref{eq:spherical_series}) of $f_\mathcal{K}$ is essentially an expansion of $f_\mathcal{K}$ in terms of the characters of irreducible representations of $G$~\cite{helgason2022groups}\cite{yaglom_second-order_1961},
\begin{equation} \label{eq:char_expansion}
    f_\mathcal{K}(x) = \sum_{\lambda \in \Lambda}\,d_\lambda\hspace{0.02cm}\hat{f}_\mathcal{K}(\lambda)\hspace{0.02cm}\chi_\lambda(x)
\end{equation}
Here, $\Lambda$ is the set of equivalence classes of irreducible representations of $G$ (not of $G \times G$), $d_\lambda$~is~the dimension of  $\lambda$ and $\chi_\lambda$ its character. The kernel $\mathcal{K}$ is universal if and only if $\hat{f}_\mathcal{K}(\lambda) > 0$ for each $\lambda \in \Lambda$. 

Recall the heat kernel of $G$ has the following expansion
\begin{equation} \label{eq:heat_lie}
    \mathcal{K}_t(x,y) = \sum_{\lambda \in \Lambda}\,d_\lambda\hspace{0.02cm}e^{-t\hspace{0.02cm}c(\lambda)}\chi_\lambda(y^{-1}x)
\end{equation}
where $c(\lambda)$ are the eigenvalues of the Laplace-Beltrami operator (\textit{i.e.} Casimir operator)~\cite{fegan}. This expansion and its various closed-form expressions have been the subject of intense study~\cite{arede}\cite{fegan}. An alternative construction of a universal kernel is given by 
\begin{equation} \label{eq:fundamental}
    f_\mathcal{K}(x) = \left|1-t\chi_f(x)\right|^{-2}
\end{equation}
where $\chi_f$ is the character of some faithful (\textit{i.e.} injective) representation $f$ of $G$ and $0< t < 1/d_f$. If~$G$ is viewed as a matrix Lie group, $\chi_f(x)$ may just be replaced with $\mathrm{tr}(x)$ (the trace of $x$). 

The point here is that any irreducible representation $\lambda$ of $G$ is contained in some tensor product of several copies of any faithful representation $f$ and its dual representation $f^*$~\cite{brocker}. Therefore, each character $\chi_\lambda(x)$ appears with a strictly positive integer weight in the decomposition of $(\chi_f(x))^n(\chi^*_f(x))^m$ (where $^*$ is the complex conjugate, and $n,m = 0, 1, 2,\ldots$) into irreducible characters. With this in mind, expand (\ref{eq:fundamental}) as a product of two binomial series, 
$$
f_\mathcal{K}(x) = \sum^\infty_{n=0}\sum^\infty_{m=0} \left((\chi_f(x))^n(\chi^*_f(x))^m\right)t^{n+m}
$$
and then decompose each product $(\chi_f(x))^n(\chi^*_f(x))^m$ into irreducible characters, to obtain
$$
f_\mathcal{K}(x) = \sum^\infty_{n=0}\sum^\infty_{m=0} \left(\sum_{\lambda \in \Lambda}c^{\lambda}_{mn} \chi_\lambda(x)\right)t^{n+m}
$$
where the sum in parentheses contains finitely many non-zero terms, and all non-zero $c^\lambda_{mn}$ are strictly positive integers~\cite{brocker} (Page 137). By rearranging this last expression, 
$$
f_\mathcal{K}(x) = \sum_{\lambda\in \Lambda} \left(\sum^\infty_{n=0}\sum^\infty_{m=0} c^{\lambda}_{mn}\hspace{0.02cm}t^{n+m}\right)\chi_\lambda(x)
$$
This is the character expansion (\ref{eq:char_expansion}) of $f_\mathcal{K}$ with all the coefficients strictly positive, as required. Therefore, the invariant kernel $\mathcal{K}(x,y) = f_\mathcal{K}(y^{-1}x)$ is indeed universal. 

\subsection{Non-compact symmetric spaces}
Theorem \ref{th:ncc0} says that the property of $C_0$-universality is shared by all invariant integrable kernels (assuming of course they are $C_0$). The task of constructing universal kernels therefore reduces to that of constructing integrable kernels. Two kinds of examples are now considered, symmetric cones and type IV symmetric spaces.

\subsubsection{Symmetric cones}
Symmetric cones are a rich subclass of the class of non-compact symmetric spaces~\cite{farautcones}. A symmetric cone is an open convex cone $X$ in some Euclidean space $V$, which has a Riemannian geometry completely determined by a so-called characteristic function
\begin{equation} \label{eq:char_func_cone}
    \Psi(x) = \int_{X}e^{-(x|y)}dy
\end{equation}
where $(x|y)$ denotes the scalar product of $x,y \in V$ and $dy$ the corresponding volume measure. In fact, $X$ becomes a Riemannian symmetric space when equipped with the metric equal to the Hessian (second derivative) of $\Psi$. This metric is given by 
\begin{equation} \label{eq:metric_cone}
    \left( u\hspace{0.02cm},v \right)_x = D^2\Psi_x(u\hspace{0.02cm},v)
\end{equation}
for $x \in X$ and $u,v \in V$ which are identified with tangent vectors to $X$ at $x$. 

For (\ref{eq:char_func_cone}) and (\ref{eq:metric_cone}) to correctly make $X$ into a Riemannian symmetric space, $X$ has to be an open convex homogeneous and self-dual cone (for the meaning of these terms, see~\cite{farautcones}). These conditions are indeed satisfied whenever $X$ is an irreducible symmetric cone, such as a~cone of real, complex, or quaternion self-adjoint positive-definite matrices of any size $r \times r$, or a Lorentz cone of any dimension $n$~\cite{farautcones} (a Lorentz cone is the cone $\Lambda_n \subset \mathbb{R}^n$ of points $x = (x_1,\ldots,x_n)$ with $x_1 > 0$ and $x^2_1 - x^2_2 - \ldots - x^2_n > 0$).

In~\cite{mostajeran2024invariantkernelsspacecomplex}, a family of integrable kernels on symmetric cones was introduced under the name of beta prime kernels. This was stated only for the cone of complex positive-definite matrices of size $r \times r$. However, both the statement and proof only relied on the symmetric cone structure.
Beta prime kernels can be described as follows.

Let $X$ be a symmetric cone as above. For $x \in X$, the ``determinant" of $x$ is $\Delta(x) > 0$, defined by
\begin{equation} \label{eq:cone_determinant}
    \Psi(x) = \Gamma_X\hspace{0.02cm}\Delta(x)^{-\frac{n}{r}}
\end{equation}
where $\Gamma_X$ is a normalising constant~\cite{farautcones} (Page 125), $n$ is the dimension of $X$ and $r$ its rank --- for a cone of positive-definite matrices, this is $r$ as above. For a Lorentz cone, it is $r = 2$.
If~$X$~is a cone of positive-definite matrices, $\Delta(x)$ is the determinant of $x$ in the usual sense. If~$X$~is a Lorentz cone $\Lambda_n \subset \mathbb{R}^{n}$, then $\Delta(x) = x^2_1 - x^2_2 - \ldots - x^2_n$ for~$x = (x_1,x_2,\ldots,x_n) \in \Lambda_n$. 

Contrary to the heat kernel, which may be quite difficult to compute, as discussed in~\cite{azangulov2023stationary}\cite{azangulov2024stationary}, the beta prime kernel has a straightforward expression,
\begin{equation} \label{eq:betaprime}
    \mathcal{K}(x,y) = \left[\frac{\Delta(x)\Delta(y)}{\Delta(x+y)^2}\right]^\alpha
\end{equation}
and is indeed $C_0$ and invariant, in addition to being integrable for $\alpha > \frac{\beta}{2}(r-1)$, where $\beta = 1, 2, 4$ when $X$ is a cone of real, complex or quaternion positive-definite matrices (respectively), and $\beta = n-2$ when $X$ is a Lorentz cone---all of these claims can be proved exactly as in~\cite{mostajeran2024invariantkernelsspacecomplex}. According to Theorem \ref{th:ncc0}, it is therefore a $C_0$-universal kernel. 

To end this discussion of symmetric cones and beta prime kernels, consider what happens in (\ref{eq:betaprime}) when $X = \Lambda_n$ is a Lorentz cone. The subset of $X$ determined by $\Delta(x) = 1$ is a hyperboloid sheet, and the Riemannian metric of $X$ (restricted to this subset) makes it a space of constant negative curvature. This will be denoted $H_{n-1} \subset \Lambda_n\hspace{0.03cm}$, since it is a hyperbolic space of dimension $n-1$.

The restriction of the integrable kernel $\mathcal{K}$ of (\ref{eq:betaprime}) to $H_{n-1}$ is again a $C_0$-universal kernel. This directly follows from the fact that it is integrable, as a kernel on $H_{n-1}\hspace{0.03cm}$, by Theorem \ref{th:ncc0}. 

For $x,y \in H_{n-1}\hspace{0.03cm}$, $\mathcal{K}(x,y)$ can be expressed in terms of the Riemannian (\textit{i.e.} hyperbolic) distance between $x$ and $y$. Specifically, note that $\Delta(x) = [x,x]$ in terms of the Minkowski form $[x,y] = x_1\hspace{0.02cm}y_1 - x_2\hspace{0.02cm}y_2 - \ldots - x_n\hspace{0.02cm}y_n\hspace{0.03cm}$. In (\ref{eq:betaprime}), $x,y \in H_{n-1}$ means $\Delta(x) = \Delta(y) = 1$. Moreover, since the Minkowski form is bilinear and symmetric,
$$
\Delta(x+y) = [x+y,x+y] = 2 + 2\hspace{0.02cm}[x,y]
$$
However, it is well-known in hyperbolic geometry that $[x,y] = \cosh(d(x,y))$ where $d(x,y)$ 
is~the~Riemannian distance~\cite{lee}. This yields the following expression 
\begin{equation} \label{eq:betaprime_hyperbolic}
    \mathcal{K}(x,y) = \frac{1}{(1+\cosh(d(x,y))^{2\alpha}}
\end{equation}
for a $C_0$-universal kernel on any hyperbolic space $H_{n-1}\hspace{0.03cm}$. When $n = 3$, which corresponds to the case of the hyperbolic plane, this is reminiscent of the Herschel-Maxwell kernel previously introduced in~\cite{dacosta2024invariantkernelsriemanniansymmetric}.

\subsubsection{Type IV symmetric spaces}
The symmetric space $X = G/K$ is called a type IV symmetric space if $G$ admits a complex Lie group structure (in other words, $G$ is the realification of a complex Lie group)~\cite{helgasonsymmetric}. 

If this is the case, then the inversion formula (\ref{eq:IFT}) lends itself to analytical evaluation, since one has the explicit expressions~\cite{helgason2022groups},
\begin{align}
\label{eq:complex_spherical}  \varphi_\lambda(\exp(a)K) &=  \frac{\Pi(\rho)}{\Pi(\mathrm{i}\lambda)}\frac{\sum_{w \in W}(\det w) e^{\mathrm{i}\lambda(w\hspace{0.02cm}a)}}{\sum_{w \in W} (\det w) e^{\rho(w\hspace{0.02cm}a)}} \\[0.1cm]
\label{eq:complex_harishchandra}
c(\lambda) & = \frac{\Pi(\rho)}{\Pi(\mathrm{i}\lambda)} 
\end{align}
where $\Pi$ is the product-of-positive-roots polynomial, given by $\Pi(\lambda) = \prod_{\alpha \in \Delta_+} (\lambda,\alpha)$ for $\lambda \in \mathfrak{a}^*_\C$
($\Delta_+ \subset \mathfrak{a}^*$ is the set of positive roots associated with the Iwasawa decomposition $G = KAN$).

Replacing (\ref{eq:complex_spherical}) and (\ref{eq:complex_harishchandra}) into (\ref{eq:IFT}) yields
\begin{equation} \label{eq:complexIFT1}
    f(\exp(a)K) = c_X\,j^{-1/2}(a)\int_{\mathfrak{a}^*}\,\hat{f}(\lambda)\left(\sum_{w \in W}(\det w) e^{\mathrm{i}\lambda(w\hspace{0.02cm}a)}\right)\Pi(-\mathrm{i}\lambda)d\lambda
\end{equation}
Here, $c_X$ is a new positive constant (different from the original one in (\ref{eq:IFT})), and $j(a)$ denotes the Weyl denominator,
\begin{equation} \label{eq:jvolf}
    j^{1/2}(a) = \sum_{w \in W} (\det w)\hspace{0.02cm} e^{\rho(w\hspace{0.02cm}a)} = \prod_{\alpha \in \Delta_+} 2\sinh(\alpha(a)/2)
\end{equation}
In dealing with (\ref{eq:complexIFT1}), it is helpful to realise this reduces to a classical inverse Fourier transform. Recall the Weyl group $W$ is a Coxeter group, generated by reflections through the hyperplanes $\alpha(a) = 0$ in $\mathfrak{a}$. The expression inside the parentheses in (\ref{eq:complexIFT1}) is skew with respect to the Weyl group (meaning a reflection applied to $a$ will change the sign of this expression). However~\cite{helgason2022groups}, the polynomial $\Pi$ is also skew with respect to the Weyl group. This implies that the ``alternating sum" can be dropped from (\ref{eq:complexIFT1}), so that
\begin{equation} \label{eq:complexIFT}
    f(\exp(a)K) = c_X\,j^{-1/2}(a)\int_{\mathfrak{a}^*}\left(\Pi(-\mathrm{i}\lambda)\hat{f}(\lambda)\right)e^{\mathrm{i}\lambda(a)}d\lambda
\end{equation}
Accordingly, the inverse spherical transform is the classical inverse Fourier transform, up to multiplication by $\Pi(-\mathrm{i}\lambda)$. This prescription can be used to produce a whole variety of invariant kernels on the symmetric space $X$. Indeed, according to the $L^1$-Godement theorem in~\cite{dacosta2024invariantkernelsriemanniansymmetric}, for~any integrable, continuous, positive-definite invariant kernel $\mathcal{K}$, the function $f_\mathcal{K}$ is given by the spherical inversion formula (\ref{eq:IFT}) (which now reduces to (\ref{eq:complexIFT})) with a positive spherical transform $\hat{f}_\mathcal{K}(\lambda)$, and any integrable, continuous, positive-definite invariant kernel is also of this form.

For example, in the case of the heat kernel $(\mathcal{P}_t\,;t > 0)$, this spherical transform is
\begin{equation} \label{eq:heatspectrum}
    \hat{f}_{\mathcal{P}_t}(\lambda) = c(t)\exp(-t((\lambda,\lambda)+(\rho,\rho)))
\end{equation}
just as in (\ref{eq:heat}). Substituting this into (\ref{eq:complexIFT}) and performing the same calculations as in~\cite{mostajeran2024invariantkernelsspacecomplex}, it~can be shown that
\begin{equation} \label{eq:complexheat}
 \mathcal{P}_t(\exp(a)K,K) = C(t)\,\frac{\Pi(a)}{j^{1/2}(a)}\exp\left[-\frac{(a,a)}{4t}\right]
\end{equation}
where $C(t)$ is a normalising constant. Since the heat kernel is invariant, this expression completely describes it in closed analytical form. The heat kernel is universal, as a direct consequence of Theorem \ref{th:ncc0}, and (\ref{eq:complexheat}) is thus a closed-form expression for a $C_0$-universal kernel. Similar expressions can be produced from (\ref{eq:complexIFT}) by replacing (\ref{eq:heatspectrum}) with other suitable choices. Intuitively, any analytic and rapidly decreasing positive $\hat{f}_\mathcal{K}(\lambda)$ will yield an integrable and therefore (by~Theorem \ref{th:ncc0}) universal kernel $\mathcal{K}$. A rigorous formulation of this intuition should be the subject of future work.

\section*{Acknowledgments}
This research is supported by the National Research Foundation, Singapore (NRF), under its NRF Fellowship Programme Award No.: NRF-NRFF18-2026-0010.

\bibliographystyle{plainnat}

\bibliography{bibliography}

\end{document}